%% file: paper.tex
\documentclass[12pt]{amsart}
\font\emailfont=cmtt10

\usepackage{amsmath,amsthm,amsfonts,amscd,flafter,epsf}

\input{basmac}

\input{macinv}

\input{macf}

\commentable{
\title[{Holomorphic triangle invariants and
symplectic four-manifolds}]
{Holomorphic triangle invariants and the topology of
symplectic four-manifolds}
\author[Peter Ozsv{\'a}th]{Peter Ozsv\'ath}
\address{Department of
Mathematics,  Princeton University, New Jersey 08540 \newline
\indent{\emailfont{petero@math.princeton.edu}}}

\author[Zolt{\'a}n Szab{\'o}]{Zolt{\'a}n Szab{\'o}} 
\address{Department of
Mathematics,  Princeton University, New Jersey 08540 \newline
\indent{\emailfont{szabo@math.princeton.edu}}}}

\begin{document}

\begin{abstract}  
	This article analyzes the interplay between symplectic
	geometry in dimension four and the invariants for smooth
	four-manifolds constructed using holomorphic triangles
	introduced in~\cite{HolDiskFour}.  Specifically, we establish
	a non-vanishing result for the invariants of symplectic
	four-manifolds, which leads to new proofs of the
	indecomposability theorem for symplectic four-manifolds and
	the symplectic Thom conjecture.  As a new application, we
	generalize the indecomposability theorem to splittings of
	four-manifolds along a certain class of three-manifolds
	obtained by plumbings of spheres. This leads to restrictions
	on the topology of Stein fillings of such three-manifolds.
\end{abstract}

\include{introsymp}

\include{toppre}
\include{adjrel}
\include{k3}
\include{nonvanish}
\include{sympthom}
\include{someplumbs}

\commentable{
\bibliographystyle{plain}
\bibliography{biblio}
}

\end{document}

%% file: basmac.tex
\newcommand\commentable[1]{#1}

\newcommand\Id{\mathrm{Id}}

\newcommand{\Tors}{\mathrm{Tors}}
\newcommand{\rk}{\mathrm{rk}}
\newcommand{\HF}{HF}

\newtheorem{theorem}{Theorem}[section]
\newtheorem{prop}[theorem]{Proposition}
\newtheorem{cor}[theorem]{Corollary}

\newtheorem{lemma}[theorem]{Lemma}

\newtheorem{defn}[theorem]{Definition}

\newtheorem{remark}[theorem]{Remark}

\def\endproof{\relax\ifmmode\expandafter\endproofmath\else
  \unskip\nobreak\hfil\penalty50\hskip.75em\hbox{}\nobreak\hfil\bull
  {\parfillskip=0pt \finalhyphendemerits=0 \bigbreak}\fi}
\def\endproofmath$${\eqno\bull$$\bigbreak}
\def\bull{\vbox{\hrule\hbox{\vrule\kern3pt\vbox{\kern6pt}\kern3pt\vrule}\hrule}}

\newcommand{\Q}{\mathbb{Q}}

\newcommand{\C}{\mathbb{C}}

\newcommand{\Z}{\mathbb{Z}}

\newcommand{\OneHalf}{\frac{1}{2}}

\newcommand{\Zmod}[1]{\Z/{#1}\Z}

\newcommand{\Ker}{\mathrm{Ker}}

\newcommand{\cm}{\cdot}

\newcommand{\Nbd}[1]{{\mathrm{nd}}(#1)}
\newcommand{\nbd}[1]{\Nbd{#1}}
\newcommand{\CDisk}{D}

\newcommand{\ModSWfour}{\mathcal{M}}
\newcommand{\ModFlow}{\ModSWfour}

\newcommand{\SpinC}{{\mathrm{Spin}}^c}

\newcommand\sgn{\mathrm{sgn}}

\newcommand\Wedge{\Lambda}

\newcommand\abuts\Rightarrow
\newcommand\Sym{\mathrm{Sym}}

\newcommand\Alg{\mathbb{A}}

%% file: macinv.tex
\newcommand\spinccanf{k}
\newcommand\spinccan{\ell}

\newcommand\HFpRed{\HFp_{\red}}
\newcommand\HFmRed{\HFm_{\red}}
\newcommand\mix{\mathrm{mix}}

\newcommand\ModSphere{\ModFlow\left({\mathbb S}\longrightarrow 
\Sym^{g-1}(\Sigma_{1})\times \Sym^2(\Sigma_{2})\right)}
\newcommand\ModSpheres\ModSphere

\newcommand\HFleq{\HF^{\leq 0}}

\newcommand\HFred{\HF_{\rm red}}

\newcommand{\red}{\mathrm{red}}

\newcommand\HFp{\HFb}
\newcommand\HFpm{HF^{\pm}}
\newcommand\HFm{\HF^-}

\newcommand\HFinf{HF^\infty}

\newcommand\HFa{\widehat{HF}}
\newcommand\HFb{HF^+}
\newcommand\gr{\mathrm{gr}}

\newcommand\UnparModSp{\widehat \ModSp}
\newcommand\UnparModFlow\UnparModSp
\newcommand\Mod\ModSp

\newcommand\PD{\mathrm{PD}}

\newcommand{\spinc}{\mathfrak s}

\newcommand{\spinct}{\mathfrak t}

\newcommand\ModMaps{\mathcal M}
\newcommand\ModSp\ModMaps

\newcommand\uHF{\underline{\HF}}

\newcommand\uHFpRed{\uHFpred}
\newcommand\uHFpred{\underline{\HF}^+_{\red}}

\newcommand\uHFinf{\uHF^\infty}

\newcommand\Fm[1]{F^{-}_{#1}}

\newcommand\Fp[1]{F^{+}_{#1}}

\newcommand\Fleq[1]{F^{\leq 0}_{#1}}

%% file: macf.tex
\newcommand{\divis}{\mathfrak d}
\newcommand\Dehn{D}
\newcommand\PiRed{\Pi^\red}

\newcommand\Fiber{F}

\newcommand\spincu{\mathfrak u}

\newcommand\Fmix[1]{F^{\mix}_{#1}}

\newcommand\Fc{F^\circ}

\newcommand\HFc{\HF^\circ}

\newcommand\Dual{\mathcal D}
\newcommand\Duality\Dual

\newcommand\Tor{\mathrm{Tor}}

%% file: introsymp.tex
\maketitle
\section{Introduction}

In~\cite{HolDiskFour}, we constructed an invariant for smooth, closed
four-manifolds (using holomorphic triangles, and the Floer
homology theories defined in~\cite{HolDisk} and \cite{HolDiskTwo}).
The aim of the present article is to investigate this invariant in the
case where $X$ is a closed, symplectic four-manifold. Our first result is the 
following:

\begin{theorem}
\label{intro:NonVanishing}
If $(X,\omega)$ is a closed, symplectic manifold with $b_2^+(X)>1$, then for
the canonical $\SpinC$ structure $\spinccanf$, we have that
$$\Phi_{X,\spinccanf}=\pm 1.$$ Moreover, if $\spinc\in\SpinC(X)$ is
any $\SpinC$ structure for which $\Phi_{X,\spinc}\not\equiv 0$, then
we have the inequality that $$\langle c_1(\spinccanf),\omega\rangle
\leq
\langle c_1(\spinc),\omega\rangle,$$
with equality iff $\spinccanf=\spinc$.
\end{theorem}

The above can be seen as a direct analogue of a theorem of Taubes
concerning the Seiberg-Witten invariants for symplectic manifolds,
see~\cite{TaubesSympI} and \cite{TaubesSympII}.  However, the proof
(given in Section~\ref{sec:NonVanishing}) is quite different in
flavor. While Taubes' theorem uses the interplay of the symplectic
form with the Seiberg-Witten equations, our approach uses the topology
of Lefschetz fibrations, together with general properties of
$\HFp$. As such, our proof relies on a celebrated result of
Donaldson~\cite{DonaldsonLefschetz}, which constructs Lefschetz
pencils on symplectic manifolds, see also~\cite{AurouxKatz} and \cite{Smith}.

Combined with the general properties of $\Phi$
(see~\cite{HolDiskFour}), the above non-vanishing theorem has a number
of consequences.

\subsection{New proofs of known results.}  

Theorem~\ref{intro:NonVanishing} can be used to reprove the
indecomposability theorem for symplectic four-manifolds, a theorem
whose K\"ahler version was established by Donaldson using his
polynomial invariants~\cite{DonaldsonPolynomials}, and whose
symplectic version was established by Taubes using Seiberg-Witten
invariants~\cite{TaubesSympI}:

\begin{cor} {\bf{(Donaldson: K\"ahler case; Taubes: symplectic case)}}
\label{intro:Indecomposability}
If $(X,\omega)$ is a closed symplectic four-manifold, then it admits no
smooth decomposition as a connected sum $X=X_1\# X_2$ into two pieces
with $b_2^+(X_1)$, $b_2^+(X_2)>0$.
\end{cor}

\begin{proof}
This follows immediately from the non-vanishing result in
Theorem~\ref{intro:NonVanishing}, together with the vanishing result
for $\Phi$ for a connected sum, 
Theorem~\ref{HolDiskFour:intro:VanishingTheorem} of~\cite{HolDiskFour}
(which in turn follows easily from the definition of $\Phi$).
\end{proof}

In the course of proving Theorem~\ref{intro:NonVanishing}, we
establish a certain ``adjunction relation'', which can be seen as an
analogue of an earlier adjunction relation from Seiberg-Witten theory
(see~\cite{FSthom} and \cite{SympThom}). Together with
Theorem~\ref{intro:NonVanishing}, this relation gives a new proof of the
symplectic Thom conjecture. Note that this question has a long history in gauge
theory. Various versions were proved in~\cite{KMPolyStruct},
\cite{KMthom}, \cite{MSzT}, and the general case (which we reprove here)
is contained in~\cite{SympThom}.

\begin{theorem} 
\label{intro:SympThom}
If $(X,\omega)$ is a symplectic four-manifold 
and $\Sigma\subset X$ is an embedded, symplectic submanifold, then
$\Sigma$ is genus-minimizing in its homology class.
\end{theorem}

\subsection{Generalized indecomposability}

We will generalize the indecomposability theorem
for symplectic four-manifolds (Corollary~\ref{intro:Indecomposability})
to a large class of plumbed three-manifolds, in place of $S^3$.

By a {\em weighted graph} we mean a graph $G$, equipped with an
integer-valued function $m$ on the vertices of $G$.
Recall
that for each weighted graph, there is a
uniquely associated three-manifold $Y(G,m)$, which is the boundary of
the associated plumbing of disk bundles over spheres (the integer
multiplicities here record the Euler numbers of the disk bundles).
The {\em degree} of a vertex $v$ in a graph $G$, denoted $d(v)$,
is the number of edges which contain 
the given vertex.

\begin{theorem}
\label{thm:CertainPlumbings}
Let $Y=Y(G,m)$ be a plumbed three-manifold, 
where $(G,m)$ satisfies the following conditions:
\begin{itemize}
\item $G$ is a disjoint union of trees
\item at each vertex in $G$, we have that
\begin{equation}
\label{eq:DegreeMultiplicity}
m(v)\geq d(v).
\end{equation}
\end{itemize} Then no closed, symplectic four-manifold $(X,\omega)$ can
be decomposed along $Y$ as a union $$X=X_1\cup_{Y} X_2$$ into two pieces with
$b_2^+(X_1)>0$ and $b_2^+(X_2)>0$.
\end{theorem}

Note that in the special cases where $Y$ is $S^2\times S^1$ or a lens space,
the above theorem was known using 
Seiberg-Witten theory.

\begin{cor}
Let $G$ be a weighted graph satisfying the hypothesis of
Theorem~\ref{thm:CertainPlumbings}. If $X$ is any Stein four-manifold
with $\partial X = \pm Y(G)$, then $b_2^+(X)=0$.
\end{cor}

\begin{proof}
According to~\cite{LMcontact}, such a Stein manifold $W$ can always be
embedded in a surface of general type $X$, so that $b_2^+(X-W)>0$.
Thus, the corollary follows from Theorem~\ref{thm:CertainPlumbings}.
\end{proof}

Note that $-Y(G)$ always admits a Stein filling with $b_2^+(X)=0$, using
a theorem of Eliashberg~\cite{EliashbergStein}, see also
\cite{GompfStipsicz}.

Theorem~\ref{thm:CertainPlumbings} follows from
Theorem~\ref{intro:NonVanishing}, coupled with a vanishing invariant
for four-manifolds admitting a decomposition along $Y(G,m)$. In turn,
this vanishing theorem follows from from a Floer homology calculation
for plumbings along graphs which satisfy the hypotheses of
Theorem~\ref{thm:CertainPlumbings}.  Of course, it is interesting to
consider plumbing diagrams which do not satisfy
Inequality~\eqref{eq:DegreeMultiplicity}. For this more general case,
one does not expect such a strong vanishing theorem -- for instance,
any Seifert fibered space with $b_1(Y)=0$ can be obtained as a
plumbing along a tree. We return to the general case of
three-manifolds obtained as plumbings along trees in a future paper,~\cite{Plumbings}.

\subsection{Organization}

This paper is organized as follows. In Section~\ref{sec:Preliminaries}
we rapidly review some of the basic notions used throughout this
paper, specifically regarding Lefschetz fibrations. We also extend the
four-manifold invariant $\Phi$ defined in~\cite{HolDiskFour} to the
case where the the four-manifold $X$ has $b_2^+(X)=1$. In
Section~\ref{sec:AdjunctionRelation}, we derive the adjunction
relation Theorem~\ref{thm:AdjunctionRelation} which is used later in
the proofs of Theorems~\ref{intro:NonVanishing} and
\ref{intro:SympThom}. In
Section~\ref{sec:K3}, we calculate $\Phi$ for the $K3$ surface. In
Section~\ref{sec:NonVanishing}, we prove
Theorem~\ref{intro:NonVanishing}, along with an auxiliary
non-vanishing result for the Floer homology groups of a three-manifold
which fibers over the circle. One ingredient in this proof is the $K3$
calculation in the previous section. In
Section~\ref{sec:SymplecticThom}, we deduce
Theorem~\ref{intro:SympThom} from Theorems~\ref{intro:NonVanishing} and
\ref{thm:AdjunctionRelation}. In Section~\ref{sec:SomePlumbings}, we provide the
Floer homology calculations which lead to  Theorem~\ref{thm:CertainPlumbings}.

This paper, of course, is built on the theory developed in~\cite{HolDisk},
\cite{HolDiskTwo}, and~\cite{HolDiskFour}, and it is written assuming familiarity
with those papers. Important properties of the four-dimensional
invariant $\Phi$ (which will be used repeatedly here) are summarized
in Section~\ref{HolDiskFour:sec:Cobordisms} of~\cite{HolDiskFour}.
Moreover, at two important points in the present paper
(when calculating the invariant for the $K3$ surface, and when 
finding examples of three-manifolds with non-trivial Floer homology
which fiber of the circle) we rely on some of the calculations of Floer homology groups
given in~\cite{AbsGraded}, see especially Section~\ref{AbsGraded:sec:SampleCalculations}
of~\cite{AbsGraded}.

\subsection{Further remarks}

For the purposes of proving Theorem~\ref{intro:SympThom}, we extend
the invariant $\Phi$ to four-manifolds with $b_2^+(X)=1$. As one
expects from the analogy with gauge theory, the invariant in that case
has additional structure. For our purposes, it suffices to construct
$\Phi$ as the invariant of a four-manifold equipped with a line $L$
inside $H_2(X;\Q)$ consisting of vectors with square zero.  This line
corresponds to a choice of a ``chamber at infinity''
(compare~\cite{Chambers}). We hope to return to this topic in a future
paper.

The pseudo-holomorphic triangles in the $g$-fold symmetric product of
the Heegaard surface implicit in the statement of
Theorem~\ref{intro:NonVanishing} naturally gives rise to a locus
inside $X$. It is quite interesting to compare this object with the
pseudo-holomorphic curve constructed by Taubes in~\cite{TaubesGeoSW}.
This may also 
provide a link with the work of Donaldson and Smith, see~\cite{DonSmith}.

\vskip.2cm
{\noindent{\bf Acknowledgements.}} It is our pleasure to thank Andr{\'a}s Stipsicz for
some very helpful discussions.

%% file: toppre.tex
\section{Preliminaries}
\label{sec:Preliminaries}

We collect here some of the preliminaries for the proof of
Theorem~\ref{intro:NonVanishing}.
In Subsection~\ref{subsec:Lefschetz}, we review some standard properties
of Lefschetz fibrations, mainly to set up the terminology which will
be used later. For a thorough discussion of this topic,
we refer the reader to~\cite{GompfStipsicz}.
We then return to some properties of $\HFpm$, building on
the results from~\cite{HolDiskFour}.

\subsection{Lefschetz fibrations}
\label{subsec:Lefschetz}

Let $C$ be an oriented two-manifold (possibly with boundary). A {\em
Lefschetz fibration over $C$} is a smooth four-manifold $W$ and a map
$\pi\colon W \longrightarrow C$ with finitely many critical points,
each of which admits an orientation-preserving chart modeled on
$(w,z)\in\C^2$, where the map $\pi$ is modeled on the map
$\C^2\longrightarrow \C$ given by $(w,z)\mapsto w^2+z^2$.
Moreover, we will always assume that any two critical points map to 
different values under $\pi$.

If $\pi\colon W \longrightarrow C$ has no critical points, then the
fibration endows $W$ with a canonical almost-complex structure,
characterized by the property that the fibers of $\pi$ are
$J$-holomorphic. Since a $\SpinC$ structure over a four-manifold is
specified by an almost-complex structure in the complement of finitely
many points, a Lefschetz fibration endows $W$ with a canonical
$\SpinC$ structure, which we denote by $\spinccanf$. We adopt here the
conventions of~\cite{TaubesSympI}: the first Chern class of the
canonical $\SpinC$ structure agrees with the first Chern class of the
complex {\em tangent} bundle (on the locus where the latter is
defined).

A Lefschetz fibration is said to be {\em relatively
minimal} if none of the fibers of $\pi$ contains exceptional spheres -- i.e.
spheres whose self-intersection number is  $-1$.

Lefschetz fibrations over the disk $\CDisk$ $$\pi\colon
W\longrightarrow \CDisk$$ (with $n$ critical points) can be specified
by an ordered $n$-tuple of simple, embedded curves $\tau_1,...,\tau_n$ in
$F$. The space $W$ then has the homotopy type of the two-complex by
attaching disks to $F$ along the curves.  Homologies between the
$[\tau_i]$ gives rise to homology classes in $W$. More precisely, we can
identify $$H_2(W;\Z)\cong \Z\oplus \Ker\left(\Z^n\longrightarrow
H_1(F;\Z)\right),$$ where the first $\Z$ factor is generated by the
homology class of the fiber $F$, and the map $\Z^n\longrightarrow
H_1(F;\Z)$ is the map generated by taking multiples of the homology
classes of $[\tau_1]$,...,$[\tau_n]$ in $H_1(F;\Z)$.

Relative minimality in this case is equivalent to the condition
that none of these distinguished curves in $F$ bound disks in $F$.

\begin{lemma}
\label{lemma:PeriodicDomainsInFibration}
Suppose that $P\subset F$ is a 
two-dimensional 
manifold-with-boundary
whose boundary is some collection
of curves among the $\{\tau_1,...,\tau_n\}$ (each with multiplicity one). Let
${\widehat P}$ denote the closed surface in $W$ obtained by
attaching copies of vanishing cycles to $P$.
Then, 
\begin{eqnarray*}
g({\widehat P})&=&g(P) \\
{\widehat P}\cm {\widehat P} &=& -\left(\#{\text{of boundary components of $P$}}\right) \\
\langle c_1(\spinccanf),[{\widehat P}] \rangle + {\widehat P}\cm{\widehat P} &=& 2-2g({\widehat P}).
\end{eqnarray*}
\end{lemma}

\begin{proof}
The equality on the genus is obvious. The self-intersection number of
${\widehat P}$ follows from the fact that the vanishing cycles are
finished off with disks with framing $-1$. The final equation is a
local calculation, in view of the fact that the determinant bundle of
the canonical $\SpinC$ structure is identified, in the complement of
the singular locus, with the bundle of fiber-wise tangent vectors.
\end{proof}

A Lefschetz fibration over a disk bounds a three-manifold which is a
surface bundle over the circle. Such a circle bundle is uniquely given
by the mapping class of its monodromy (a mapping class of a
two-manifold is an orientation-preserving diffeomorphism, modulo
isotopy). 
Recall that  a {\em
(right-handed) Dehn twist} of the annulus (using the conventions
of~\cite{GompfStipsicz}) is a diffeomorphism $\Psi$ of $[0,1]\times S^1$
which fixes the boundary pointwise, and satisfies the additional
property that the intersection number of an
$$\#\Big([0,1]\times\{x\}\cap
\psi([0,1]\times\{x\})\Big)=-1.$$ More generally, a (right-handed)
Dehn twist about a curve $\tau\subset F$ is a self-diffeomorphism
$\Dehn_{\tau}$ of $F$ whose restriction to some annular neighborhood of
$\tau$ is a right-handed Dehn twist of the annulus, and which fixes
all points in the complement in $F$ of the annular neighborhood. If
the Lefschetz fibration has a unique critical point, then its
monodromy is a Dehn twist about some curve $\tau$ in the fiber $F$.
More generally, if the fibration has critical values
$\{x_1,...,x_n\}$, then we can find the tuple of curves
$(\tau_1,...,\tau_n)$ by embedding a bouquet of $n$ circles in
$\CDisk-\{x_1,...,x_n\}$, so that the winding number of $\tau_i$
around $x_j$ is $\delta_{i,j}$. Then, the monodromy about the
$i^{th}$ circle is a Dehn twist about $\tau_i$. Thus, the monodromy 
map around the boundary of the disk is given as the product of Dehn
twists $\Dehn_{\tau_1}\circ...\circ\Dehn_{\tau_n}$.

Note that the curves $(\tau_1,...,\tau_n)$
obtained from a Lefschetz fibration as above depend on the embedding
of the bouquet of circles. By changing the homotopy classes of the
embedded circles, we can vary the curves $(\tau_1,...,\tau_n)$ by {\em
Hurwitz moves}, moves which carry the tuple
$(\tau_1,...,\tau_i,\tau_{i+1},...,\tau_n)$ to
$(\tau_1,...,\tau_{i+1},
\Dehn_{\tau_{i+1}}(\tau_{i}),...,\tau_n)$. 

It is well-known that any orientation-preserving
automorphism of $F$ extends to a Lefschetz fibration over the disk.
Indeed, we find it convenient to formulate this fact as follows:

\begin{theorem} (see~\cite{Humphries})
\label{thm:GenMappingClassGroup}
The mapping class group is generated as a monoid by Dehn twists along
finitely many  non-separating
curves. Indeed we can choose the generating set
$\{\tau_1,...,\tau_m\}$ so that their homology classes span $H_1(\Sigma;\Z)$,
all homological relations between the
curves are generated (over $\Z$) by special relations in which the
homology classes of $[\tau]_i$ appear with multiplicities zero or $1$,
and
the curves which appear with-non-zero multiplicities in these relations can be chosen to
be disjoint from one another.
\end{theorem}

\begin{proof}
It is a theorem of Humphries (see~\cite{Humphries}) that the mapping
class group is generated (as a group) by the $2g+1$ curves
$\{\alpha_1,...,\alpha_{g},\beta_1,...,\beta_g,\delta\}$ which are
pictured in Figure~\ref{fig:MappingClassGenerators}. 


\begin{figure}
\mbox{\vbox{\epsfbox{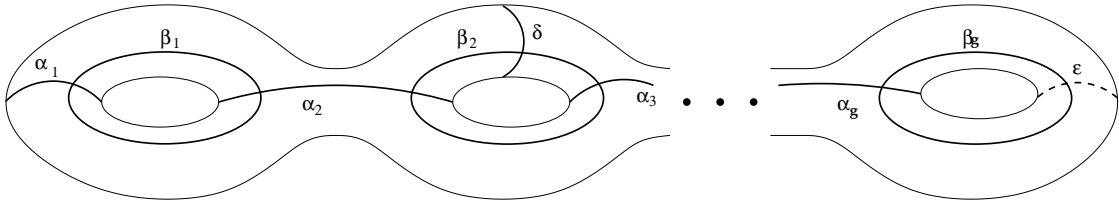}}}
\caption{\label{fig:MappingClassGenerators}
{\bf{Generators of the mapping class group.}}  Dehn twists about the
pictured curves $\{\alpha_1,...,\alpha_g,\beta_1,...,\beta_g,\delta\}$
generate the mapping class group. The additional curve $\epsilon$ is
discussed in the proof of Theorem~\ref{thm:GenMappingClassGroup}.}
\end{figure}

Now, it is easy
to see that if we include in addition the curve $\epsilon$, then we
can express the inverses of Dehn twists along all of the $\alpha_i$
and $\beta_j$ as positive multiples of Dehn twists along copies of all
the $\alpha_i$, $\beta_j$, and $\epsilon$. 
This can be seen, for example, from the
identity: : $$1 = \left(\left(\prod_{i=1}^g \Dehn_{\alpha_i}\cm
\Dehn_{\beta_i}\right)\cm \Dehn_{\epsilon}^2 \cm
\left(\prod_{i=1}^g \Dehn_{\beta_{g-i+1}}\cm \Dehn_{\alpha_{g-i+1}}\right)\right)^{4},
$$ which in turn can be obtained by exhibiting a Lefschetz fibration
over the two-sphere whose monodromy representation is given by the
above curves.  (That Lefschetz fibration is obtained by viewing the
elliptic surface $E(2g)$ as a genus $2g$ fibration over the two-sphere
-- see Chapter~8 of~\cite{GompfStipsicz} for an extensive discussion).
It remains to capture $\delta^{-1}$. To this end, we observe that $F$
has a rotational symmetry $\phi\colon F\longrightarrow F$ with the
property that we can introduce a new curve $\alpha_{g+1}$ 
so that for $i=1,...,g$, 
$\Psi(\beta_i)=\beta_{j}$ where
$j\equiv i+1\pmod{g}$, 
for $i=2,...,g$, $\Psi(\alpha_i)=\alpha_{i+1}$,
$\Psi(\alpha_{g+1})=\alpha_2$, 
$\Psi(\epsilon)=\alpha_1$, and finally
$\Psi(\alpha_{1})=\delta$. It is now clear that the mapping class
group is generated as a monoid by Dehn twists about the $2g+3$ curves
$\{\alpha_1,...,\alpha_{g+1},\beta_1,...,\beta_g,\delta,\epsilon\}$.
For homological relations between these curves, observe that the
homology classes of the
$\{\alpha_1,...,\alpha_g,\beta_1,...,\beta_g\}$ span
$H_1(\Sigma;\Z)$. It follows that the following three relations span 
all relations:
\begin{eqnarray*}
{[\alpha_1]}+[\alpha_2]+[\delta]&=&0, \\
{[\epsilon]}+[\alpha_{g+1}]+[\alpha_1]&=&0, \\
{[\alpha_2]}+...+[\alpha_{g+1}]&=& 0
\end{eqnarray*}
(See Figure~\ref{fig:MonoidGen} for an illustration in the case where $g=4$.)
\begin{figure}

\mbox{\vbox{\epsfbox{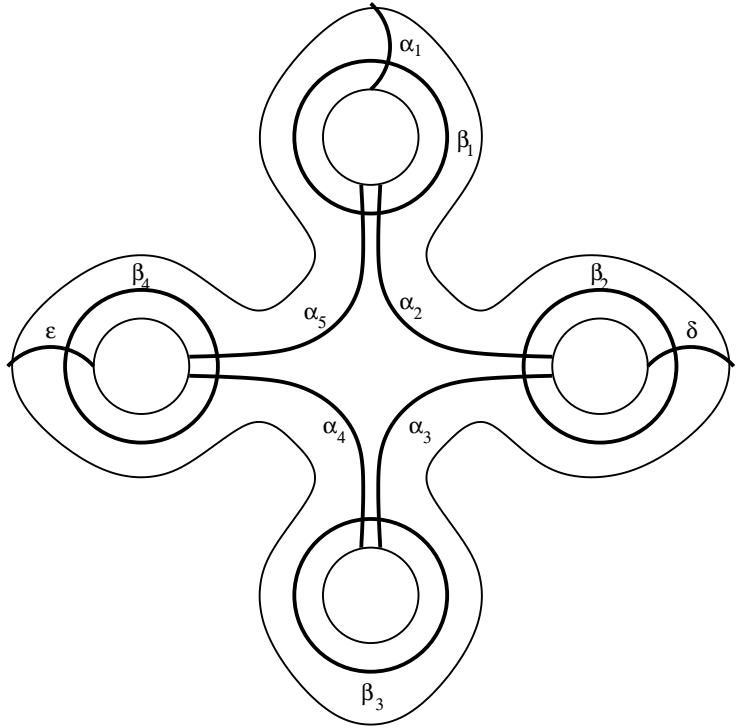}}}
\caption{\label{fig:MonoidGen}
{\bf{Monoid generators of the mapping class group, $g=4$.}}
Dehn twists about the pictured curves 
$\{\alpha_1,...,\alpha_5,\beta_1,...,\beta_4,\epsilon,\delta\}$ 
generate the mapping class group as a monoid. The symmetry $\Psi$ 
described in the proof of Theorem~\ref{thm:GenMappingClassGroup}
is realized by a $90^{\circ}$ clockwise rotation of this picture.}
\end{figure}
\end{proof}

Recall that a $\SpinC$ structure over a three-manifold $Y$ is a
suitable equivalence class of nowhere vanishing vector field over
$Y$. A three-manifold which fibers over the circle has a canonical
$\SpinC$ structure, induced by a vector field which is everywhere
transverse to the fibers. When $Y$ bounds a Lefschetz fibration over a
disk, this $\SpinC$ structure is the restriction of the canonical
$\SpinC$ structure of the Lefschetz fibration.

\subsection{Symplectic manifolds and Lefschetz fibrations}

A symplectic structure on a four-manifold $(X,\omega)$ gives the
manifold an isotopy class of
almost-complex structures, and hence a canonical
$\SpinC$ structure. Symplectic manifolds can be blown up, to construct
a new four-manifold ${\widehat X}$, which is diffeomorphic to the
connected sum of $X$ with the complex projective plane given the
opposite of its complex orientation. Symplectically, ${\widehat X}$ is
obtained by gluing the complement of a ball in $X$ to a neighborhood
of a symplectic two-sphere $E$ with self-intersection number
$-1$. Note that the canonical $\SpinC$ structure ${\widehat
\spinccanf}$ is the $\SpinC$ structure which agrees with $\spinccanf$
in the complement of $E$, and which satisfies
$$\langle c_1({\widehat \spinccanf}), [E]\rangle = +1.$$

In~\cite{DonaldsonLefschetz}, Donaldson showed that if $(X,\omega)$ is
a symplectic four-manifold, then after blowing up $X$ sufficiently
many times, one obtains a new symplectic four-manifold $({\widehat
X},{\widehat \omega})$ which admits a Lefschetz fibration $$ \pi\colon
{\widehat X}\longrightarrow S^2. $$ In fact, the fibers of $\pi$ are
symplectic, and hence the canonical $\SpinC$ structure of the
symplectic form agrees with the canonical class of the Lefschetz
fibration in the sense of Subsection~\ref{subsec:Lefschetz}.

\subsection{Preliminaries on $\HFp$}

Let $\spinct$ be a $\SpinC$ structure on an oriented three-manifold
$Y$. If $c_1(\spinct)$ is a torsion class, we simply call $\spinct$ a
torsion homology class.  The {\em divisibility} of a $\spinc$
structure is the quantity defined by 
$$\divis(\spinct)=\gcd_{\xi\in H^1(Y;\Z)} \langle c_1(\spinct)\cup
\xi, [Y]\rangle.$$

\begin{lemma}
\label{lemma:AnnihilateHFinf}
Let $Y$ be a three-manifold equipped with a non-torsion $\SpinC$ structure
$\spinct$, and let $\divis(\spinct)=d$ denote its divisibility,
then
$$(1-U^{d/2})\HFinf(Y,\spinct)=0.$$
\end{lemma}

\begin{proof}
This is an easy consequence of the material in
Section~\ref{HolDiskTwo:sec:HFinfty}
of~\cite{HolDiskTwo}. Specifically, it is shown there
(Theorem~\ref{HolDiskTwo:thm:HFinfTwist}) that the twisted version of
$\HFinf$ is, 
$\uHFinf(Y,\spinct)$
is a free $\Z[U,U^{-1}]$ module, 
endowed with the $\Z[H^1(Y;\Z)]$ action where $e^h$ ($h\in H^1(Y;\Z)$)
acts as multiplication by $U^{\langle h\cup c_1(\spinct),
[Y] \rangle/2}$. There is a universal coefficients spectral sequence converging
to the untwisted version $\HFinf(Y)$ (as a $\Z[U,U^{-1}]$ module), 
and whose $E_2$ term is given by
$$\Tor^i_{\Z[U,U^{-1}]}(\uHFinf_j(Y,\spinct),\Z[U,U^{-1}]),$$
where here the $\Z[U,U^{-1}]$ is given a trivial action by $\Z[H^1(Y;\Z)]$.
Observe that we have a free resolution of
$\uHFinf_j(Y,\spinct)$ as a module over $\Alg=\Z[U,U^{-1}]\otimes_{\Z}\Z[H^1(Y;\Z)]$, given by
$$
\begin{CD}
\bigotimes_{i=1}^{b_1(Y)} \Alg @>{e^{h_i}-U^{n_i/2}}>>\Alg,
\end{CD}
$$
where $h_i$ is a basis for $H^1(Y;\Z)$, and $n_i=\langle c_1(\spinct)\cup h_i,
[Y]\rangle$. So, the $E_2$ term of the above sequence is simply calculated by
the homology of
$$
\left( \begin{CD}
\bigotimes_{i=1}^{b_1(Y)} 
\Z[U,U^{-1}] @>{1-U^{n_i/2}}>> \Z[U,U^{-1}].
\end{CD}
\right).
$$
Bearing in mind that
$$\left(\frac{\Z[U]}{U^{a}-1}\right)
\otimes_{\Z[U,U^{-1}]}
\left(\frac{\Z[U]}{U^{b}-1}\right)
\cong \Z[U]/(U^{c}-1) 
\cong  \Tor^1_{\Z[U,U^{-1}]}\left(\frac{\Z[U]}{U^{a}-1},
\frac{\Z[U]}{U^{b}-1}\right)$$
(and all higher $\Tor^i$ vanish), where here $c=\gcd(a,b)$, 
it follows easily that
$U^{d/2}-1$ annihilates this $E_2$ term (in view
of the fact that $d$ is the greatest common divisor of
the integers $\langle c_1(\spinct)\cup
h_i,[Y]\rangle/2$ for $i=1,...,b_1(Y)$), 
and hence it also annihilates $\HFinf(Y)$ with untwisted coefficients.
\end{proof}

Let $Y$ be a closed, oriented three-manifold. It follows from
Lemma~\ref{lemma:AnnihilateHFinf} and the finiteness of $\HFred(Y)$
that for any sufficiently large integer $k$ so that if $\spinct$ is a
non-torsion $\SpinC$ structure with divisibility $d$, then
$$(1-U^{dk/2})\colon
\HFm(Y,\spinct) \longrightarrow \HFmRed(Y,\spinct)$$ defines a
projection map of $\HFm(Y,\spinct)$ onto $\HFmRed(Y,\spinct)$. 
In fact, by composing with the inverse of the coboundary map
$$\tau\colon \HFpRed(Y,\spinct)\longrightarrow \HFmRed(Y,\spinct),$$
this gives a map
$$\PiRed_Y\colon \HFm(Y,\spinct)\longrightarrow \HFpRed(Y,\spinct).$$

Using a decomposition of $W$ along such a three-manifold $N$ (and
using a $\SpinC$ structure $\spinc$ over $W$ whose restriction to $N$
is non-torsion) is analogous to the ``admissible cuts''
of~\cite{HolDiskFour}. Indeed, the comparison with the mixed
invariants defined there is given by the following:

\begin{prop}
\label{prop:OtherCut}
Suppose that $W$ is a cobordism from $Y_1$ to $Y_2$ with $b_2^+(W)>1$,
which is separated by a three-manifold $N$ into a pair of cobordisms
$W_1\cup_N W_2$. Given any pair of $\SpinC$ structures $\spinc_1$ and $\spinc_2$ 
over $W_1$ and $W_2$ respectively whose restrictions to $N$ agree and are non-torsion,
we have:
$$\Fp{W_2,\spinc|W_2}\circ \PiRed_N\circ
\Fm{W_1,\spinc|W_1}(\xi)=
\sum_{\{\spinc\in\SpinC(W)|\spinc|W_1=\spinc_1,\spinc|W_2=\spinc_2\}}
\pm \Fmix{W,\spinc}(\xi).$$
\end{prop}

\begin{proof}
Since $c_1(\spinc)|N$ is non-torsion, we can find an embedded surface
$F\subset N$ with $\langle c_1(\spinc),[F]\rangle\neq 0$.  Now, we can
cut $W$ in two along $N'=Y_1\# (S^1\times F)$, giving $W=W_1'\cup_{N'}
W_2'$.  Now, by naturality of the exact sequences (relating $\HFm$,
$\HFinf$, and $\HFp$) the usual composition laws, we see that
$$\Fp{W_2,\spinc|W_2}\circ \PiRed_N\circ
\Fm{W_1,\spinc|W_1}(\xi)=
\sum_{\eta\in \delta H^1(N)}\Fp{W_2',\spinc+\eta|W_2'}\circ 
\PiRed_{N'}\circ \Fp{W_1',\spinc+\eta|W_1'}(\xi).$$

Next, we find some embedded surface $\Sigma\subset W$ of positive
square which  disjoint from $F$, and let $Q$ denote its tubular neighborhood.
Then, $Q\# Y_2$ naturally gives a cut of 
$W$ which we can arrange to be
disjoint
from the cut $N'$ used above (by making the tubular neighborhoods 
sufficiently small). It then follows now easily from
from the composition laws that
$$
\sum_{\eta\in \delta H^1(N)}\Fp{W_2',\spinc+\eta|W_2'}\circ 
\PiRed_{N'}\circ \Fp{W_1',\spinc+\eta|W_1'}(\xi)
=\sum_{\eta\in\delta H^1(N)}
\Fmix{W,\spinc+\delta\eta}\left((1-U^{dk/2}\right)\xi).  $$
The equation follows by choosing $k$ large enough that $U^{dk/2}$
annihilates all the mixed invariants of $W$.
\end{proof}

\subsection{The case where $b_2^+(X)=1$}
\label{subsec:BTwoOne}

The construction of closed invariants defined in~\cite{HolDiskFour}
works only in the case where the four-manifold has
$b_2^+(X)>1$. However, Proposition~\ref{prop:OtherCut} suggests a
construction which can be used even when $b_2^+(X)=1$. Rather than
setting up the general theory at present, we content ourselves with
developing enough of it to allow us to establish
Theorem~\ref{intro:SympThom} in the case where $b_2^+(X)=1$.

\begin{defn}
\label{def:Cuts}
Let $X$ be a closed, smooth four-manifold and choose a line
$L\subset H_2(X;\Q)$ with the property that each vector $v\in L$
has $v\cm v=0$. Choose a cut $X=X_1\#_N X_2$ for which 
the image of $H_2(N;\Q)$ inside $H_2(X;\Q)$ is $L$. 
Then, for each $\SpinC$ structure
$\spinc\in\SpinC(X)$ for which $c_1(\spinc)$ evaluates non-trivially
on $L$, we can define 
$$\Phi_{W,\spinc,L}\colon \Z[U]\otimes \Wedge^*(H_1(X)/\Tors) 
\longrightarrow \Z/\pm 1$$
to be the non-zero on only those homogeneous elements of
$\Z[U]\otimes \Wedge^*(H_1(X)/\Tors)$ whose
degree is given by
$$d(\spinc)=\frac{c_1(\spinc)^2-2\chi(X)-3\sgn(X)}{4},$$
where $\chi(X)$ denotes the Euler characteristic of $X$
and $\sgn(X)$ denotes the signature of its intersection form.
On those elements, the invariant is
the coefficient of $\Theta^+\in\HFp(S^3,\spinc)$
in the expression
$$\Fp{W_2,\spinc|W_2}\circ \PiRed_N\circ
\Fm{W_1,\spinc|W_1}(U^n\cm \Theta^-\otimes \zeta).$$
Here, $\Theta^+$ and $\Theta^-$ are bottom- and top-dimensional
generators of $\HFp(S^3)$ and $\HFm(S^3)$ respectively.
\end{defn}

\begin{prop}
The invariant $\Phi_{W,\spinc,L}$ depends on the cut only
through the choice of line $L\in H_2(X;\Q)$.
\end{prop}

\begin{proof}
An embedded surface $F\subset X$ whose homology class is in the line $L$
always gives rise to a cut as in Definition~\ref{def:Cuts}. Specifically,
let $F\subset X$ be a smoothly embedded, connected submanifold with $[F]\in L$.
Then, we decompose 
$$
X=\left(X-\nbd{F}\right)\cup_{S^1\times F}\left(F\times\CDisk\right).$$

Next, suppose that $F_1$ and $F_2$ are two embedded surfaces whose
homology classes lie inside $L$.  Then we claim that there is a
third embedded surface $F_3$ which is disjoint from both $F_1$ and
$F_2$, and whose homology class also lies inside $L$. This is
easily constructed by starting with some initial surface $\Sigma$, and
then adding handles along canceling pairs of intersection points
between $\Sigma$ and $F_1$ (and then $\Sigma$ and $F_2$). It follows now from the 
usual arguments that the invariant calculated by using 
the cut determined by $F_1$ (or $F_2$) agrees with the invariant calculated 
using the cut determined by $F_3$; i.e. the invariant using any such embedded
surface is independent of the choice of homology class and surface.

Finally, if $X=W_1\cup_N W_2$ is an arbitrary cut as in Definition~\ref{def:Cuts},
then we can find an embedded surface $F\subset X$ disjoint from $N$ whose
homology class lies in the line $L$. Indeed, 
letting $F_0$ be any surface representing an
element of $H_2(N;\Z)$ with non-trivial image in $H_2(X;\Z)$, we let $F$ be a
surface obtained by pushing $F_0$ out of $N$, using some vector field normal to $N$
inside $X$. Since $F$ is disjoint from $N$, again, the usual arguments show
that the invariant calculated using the cut $N$ agree with the invariants
calculated using the cut determined by any embedded surface whose homology class
lies in $L$.
\end{proof}

%% file: adjrel.tex
\section{The adjunction relation}
\label{sec:AdjunctionRelation}

We prove here  the following adjunction relation
(for the Seiberg-Witten analogues,
compare~\cite{FSthom} when $g=0$, and \cite{SympThom} when $g>0$):

\begin{theorem}
\label{thm:AdjunctionRelation}
For each genus $g$ there is an element $\xi\in\Z[U]\otimes_\Z \Wedge^*
H_1(\Sigma)$ of degree $2g$ with the following significance. Given any
smooth, oriented, four-dimensional cobordism $W$ from $Y_1$ to $Y_2$
(both of which are connected three-manifolds), any smoothly-embedded
connected, oriented submanifold $\Sigma\subset W$ of genus $g$, and
any $\spinc\in\SpinC(W)$ satisfying the constraint that 
\begin{equation}
\label{eq:AdjIneqViolated}
\langle
c_1(\spinc),[\Sigma]\rangle-[\Sigma]\cm[\Sigma] = -2g(\Sigma),
\end{equation}
then we have the
relation: 
\begin{equation}
\label{eq:AdjRel}
\Fc_{W,\spinc}(\cdot)=
\Fc_{W,\spinc+\epsilon \PD[\Sigma]}(i_*(\xi(\Sigma))\otimes \cdot),
\end{equation}
where $\epsilon$ is the sign of $\langle c_1(\spinc),[\Sigma]\rangle$,
and $i_*\colon \Z[U]\otimes_\Z \Wedge^* H_1(\Sigma) \longrightarrow
\Z[U]\otimes_\Z \Wedge^* H_1(W)/\Tors$ is the map 
induced by the inclusion $i\colon \Sigma\longrightarrow W$.
\end{theorem}

Before proceeding to the proof of
Theorem~\ref{thm:AdjunctionRelation}, we make a few general
observations. Note that if
$$\Big|\langle c_1(\spinc),[\Sigma]\rangle\Big| \geq 2g-\Sigma\cm\Sigma,$$
the above theorem always obtain relations of the form 
of Equation~\eqref{eq:AdjRel}, which can be obtained 
by reversing the
orientation of $\Sigma$ and adding extra 
null-homologous handles if necessary, to achieve the hypotheses
of Theorem~\ref{thm:AdjunctionRelation}.

It is not important for our present purposes to identify the
particular word $\xi(\Sigma)$.
However, it is easy to see that for a genus $g$ surface,
$$\xi(\Sigma)\equiv U^{g} \pmod{\Wedge^* H_1(\Sigma)},$$ 
by observing that surfaces and $\SpinC$ structures
satisfying the hypotheses of Theorem~\ref{thm:AdjunctionRelation}
can be found
in a tubular neighborhood of a two-sphere of arbitrary negative
self-intersection number, where all the maps on $\HFinf$ are
non-trivial. Indeed, it is natural to expect from the
analogy with Seiberg-Witten theory that $\xi(\Sigma)$
is given by the 
formula:
$$\xi(\Sigma)=\prod_{i=1}^{g}(U-A_i\cm B_i)$$ (compare~\cite{SympThom}).

With these remarks in place, we turn our attention to the proof
of Theorem~\ref{thm:AdjunctionRelation}.
One ingredient in this proof is the behavior of $\HFc$ under connected
sums, as we recall presently.  In Section~\ref{AbsGraded:sec:CorrTerm}
of~\cite{AbsGraded}, we defined a product $$\bigotimes\colon
\HFc(Y,\spinct)\otimes_{\Z[U]} \HFleq(Z,\spincu)
\longrightarrow \HFc(Y\#Z, \spinct\#\spincu),$$
which, in the case where $Z\cong S^3$ is an isomorphism (indeed,
it is the canonical isomorphism obtained from the diffeomorphism
$Y\# S^3$ with $Y$).
This
product is functorial under cobordisms (see
Proposition~\ref{AbsGraded:prop:ConnSumTransformation}
of~\cite{AbsGraded}), in the sense that if $W$ is a cobordism from
$Z_1$ to $Z_2$ equipped with the $\SpinC$ structure $\spinc$, then the
following diagram commutes: 
\begin{equation}
\label{eq:KunnethNat}
\begin{CD}
\HFc(Y)\otimes_{\Z[U]}\HFleq(Z_1) @>{\otimes}>>
\HFc(Y\# Z_1,\spinct\#\spincu_1) \\
@V{\Id\otimes \Fleq{W,\spinc}}VV
@VV{\Fc_{([0,1]\times Y)\# W,\spinct\#\spincu}}V \\
\HFc(Y)\otimes_{\Z[U]}\HFleq(Z_2) @>{\otimes}>>
\HFc(Y\# Z_2,\spinct\#\spincu_2).
\end{CD}
\end{equation}
In the above diagram, $([0,1]\times Y)\# W$ denotes the boundary connected sum.

\vskip.2cm
\noindent{\bf{Proof of Theorem~\ref{thm:AdjunctionRelation}.}}
By the blowup formula, it suffices to consider the case where
$$\Sigma\cm \Sigma =-n,$$
where $n\geq 2g$.

Now, let $N$ be a tubular neighborhood of an oriented two-manifold of
genus $g$ with self-intersection number $-n\leq -2g$, and let
${\mathfrak u}$ denote the $\SpinC$ structure over $N$ with $$\langle
c_1(\spincu),[\Sigma]\rangle = -n-2g$$ An easy application of the long
exact sequence for integral surgeries, together with the adjunction
inequality for three-manifolds (see
Theorems~\ref{HolDiskTwo:thm:ExactP} and
\ref{HolDiskTwo:thm:Adjunction} of~\cite{HolDiskTwo} respectively),
gives us that $$\HFp(Z,\spincu|Z)\cong
\Z[U^{-1}]\otimes \Wedge^*H^1(\Sigma_g).$$
(Details are given in 
Lemma~\ref{AbsGraded:lemma:CorrTermCircleBundle}
of~\cite{AbsGraded},
where the absolute grading on $\HFc(Z,\spincu|Z)$
is also calculated.) 
In particular, 
$\HFpRed(Z,\spincu|Z)=0$, and hence
$$\HFleq(Z,\spincu|Z)\cong 
\Z[U]\otimes \Wedge^*H^1(\Sigma_g).$$
Indeed, since 
$$\langle c_1(\spincu-\PD[\Sigma])^2 ,[N]\rangle 
> \langle c_1(\spinc')^2 ,[N]\rangle $$
for any $\spinc'\in\SpinC(N)$ with $\spinc'\neq\spincu - \PD[\Sigma]$
and $\spinc'|Z=\spincu|Z$, we have that
the map
$$\Fleq{N,\spincu-\PD[\Sigma]}\colon \HFleq(S^3) \longrightarrow 
\HFleq(Z,\spincu)$$
takes a top-dimensional $\Theta_{S^3}$ of $\HFleq(S^3)$ to a
top-dimensional generator $\Theta_{Z}$ of
$\HFleq(Z,\spincu|Z)$.  Moreover, according to the dimension formula,
the grading of $\Fleq{N,\spincu}(\Theta_{S^3})$ is $2g$ less than
the grading of this element so (since $\HFleq(Z,\spincu|Z)$ is
generated by $\Theta_Z$ as a module over the ring
$\Z[U]\otimes \Wedge^* H_1(Z)/\Tors\cong \Z[U]\otimes \Wedge^* H_1(\Sigma)$) 
we can find an element $\xi(\Sigma)$ of degree $2g$
in the graded algebra
$\xi(\Sigma)\in \Z[U]\otimes_\Z \Wedge^* H_1(\Sigma)$ with the property that
$$\Fleq{N,\spincu}(\Theta_{S^3})=
\xi(\Sigma)\cm \Fleq{N,\spincu-\PD[\Sigma]}(\Theta_{S^3}).$$

Next, suppose that  $Y_1$ is a three-manifold equipped with the $\SpinC$ structure
$\spinct_1$, and $W_1$ is the connected sum $\left([0,1]\times Y_1\right)\# N$,
then the naturality of the product map
(Diagram~\eqref{eq:KunnethNat}) shows that
\begin{eqnarray*}
\Fc_{W_1,\spincu}(\zeta)&=&\zeta\otimes
\Fleq{N,\spincu}(\Theta_{S^3}) \\ &=& \zeta \otimes \left(\xi(\Sigma)\cm
\Fleq{N,\spincu-\PD[\Sigma]}(\Theta_{S^3})\right) \\ &=&
\Fc_{W_1,\spincu-\PD[\Sigma]}(\xi(\Sigma)\otimes \zeta).
\end{eqnarray*}

Finally, if $W$ is a cobordism as in the statement of the theorem, we
can decompose it into a union of $W_1$ (the connected sum of a collar
neighborhood of $Y_1$ with a tubular neighborhood $N$ of $\Sigma$) and
its complement $W_2$. Both $\SpinC$ structures $\spinc$ and
$\spinc-\PD[\Sigma]$ agree over $W_2$, so the theorem follows from the
above equation, together with the composition law for the cobordism
invariant.
\qed
\vskip.2cm

%% file: k3.tex
\section{The invariant for the $K3$ surface}
\label{sec:K3}

In proving the non-vanishing theorem for symplectic four-manifolds in
general, it is helpful to have one explicit example. The aim of the
present subsection is such a calculation, for the $K3$ surface.
Recall that the $K3$ surface is the simply-connected smooth
four-manifold which can be given the structure of a compact algebraic
surface whose canonical class is trivial -- i.e. if $\spinccanf$ is the
canonical $\SpinC$ structure coming from the almost-complex structure,
then $c_1(\spinccanf)=0$.

\begin{prop}
\label{prop:K3Surface}
The invariants for the $K3$ surface are given by:
$$\Phi_{K3,\spinc}=\left\{\begin{array}{cc}
1 & {\text{if $c_1(\spinc)=0$}} \\
0 & {\text{otherwise.}}
\end{array}\right.$$
\end{prop}

We model this calculation on a paper by Fintushel and Stern
(see~\cite{FintSternK3}) where they calculate a Donaldson invariant
for $K3$ using Floer's exact triangle.  In particular, they employ the
following handle decomposition of $K3$.

Following the notation of~\cite{FintSternK3}, let $M\{p,q,r\}$ denote
the three-manifold obtained by surgeries on the Borromean rings, with integer 
coefficients $p$, $q$ and $r$. There
is a cobordism $X$ from the Poincar\'e homology three-sphere
$\Sigma(2,3,5)\cong M\{-1,-1,-1\}$ to itself with the opposite
orientation, $-\Sigma(2,3,5)=M\{1,1,1\}$, 
composed of six two-handles apiece, which we break up as 
the following composition: $$
M\{-1,-1,-1\}\Rightarrow M\{-1,-1,0\}\Rightarrow
M\{-1,-1,1\}\Rightarrow $$ $$ M\{-1,0,1\}\Rightarrow
M\{-1,1,1\}\Rightarrow M\{0,1,1\}\Rightarrow M\{1,1,1\}.$$ The
two-handles are attached in the obvious manner: for example, to go
from $M\{p,q,r\}$ to $M\{p+1,q,r\}$, we attach a two-handle along an
unknot with framing $-1$ which links the first ring once.
Let $E$ denote the negative-definite manifold obtained
as a plumbing of two-spheres according to the $E8$ Dynkin diagram;
then $\partial E = M\{-1,-1,-1\}$. There is a decomposition of $K3$ as
$$K3 \cong E\# X \# E.$$

To obtain an admissible cut of the $K3$ as required in the definition
of $\Phi$ (c.f.  Definition~\ref{HolDiskFour:def:AdmissibleCut}
of~\cite{HolDiskFour}), we cut the surface
along $N=M\{-1,-1,1\}$, to get the 
decomposition of $K3-B^4-B^4$ as
$$X_1 = \Big(S^3\Rightarrow M\{-1,-1,-1\}\Rightarrow M\{-1,-1,0\}\Rightarrow
M\{-1,-1,1\}\Big).$$
and 
$$X_2 = \Big(M\{-1,-1,1\}\Rightarrow 
M\{-1,0,1\}\Rightarrow M\{-1,1,1\}\Rightarrow M\{0,1,1\}\Rightarrow
M\{1,1,1\}\Rightarrow S^3\Big).$$

Our goal now is to determine the maps on Floer homology induced by
these two-handle additions. Indeed, the Floer homology groups
themselves, as absolutely graded groups, were calculated in
Section~\ref{AbsGraded:sec:SampleCalculations} of~\cite{AbsGraded}. In particular,
it is shown there that
\begin{eqnarray*}
\HFp_k(M\{1,1,1\})&\cong&
\left\{\begin{array}{ll}
\Z	&	{\text{if $k$ is even and $k\geq 2$}} \\
0	&	{\text{otherwise}},
\end{array}
\right.
\\
\HFp_k(M\{0,1,1\})&\cong&
\left\{\begin{array}{ll}
\Z	&	{\text{if $k\equiv \OneHalf\pmod{2}$ and
		$k \geq {\frac{3}{2}}$}} \\
0	&	{\text{otherwise}},
\end{array}
\right.
\\
\HFp_k(M\{-1,1,1\})&\cong&
\left\{\begin{array}{ll}
\Z\oplus \Z & 	{\text{if $k=0$}} \\
\Z	&	{\text{if $k$ is even and $k>0$}} \\
0	&	{\text{otherwise}},
\end{array}
\right. \\
\HFp_k(M\{-1,0,1\})
&\cong & \left\{\begin{array}{ll}
\Z & {\text{if $k\equiv \OneHalf\pmod{\Z}$ and $k\geq \OneHalf$}} \\
\Z\oplus \Z &{\text{if $k=-\OneHalf$}} \\
0 & {\text{otherwise.}}
\end{array}
\right.
\end{eqnarray*}
It is also shown there that the $\Z[U]$ action is surjective for the
first two examples, while it has a one-dimensional cokernel for the
second two.  The groups $\HFm_k$ for these three-manifolds can be
immediately deduced by the long exact sequence relating $\HFm$ and
$\HFp$ (see~\cite{HolDiskTwo}), and the groups for the remaining
three-manifolds are determined by the duality of $\HFpm$ under
orientation reversal, and the observation that
$-M\{p,q,r\}\cong M\{-p,-q,-r\}$.

In the above statements, we are using the absolute gradings on the
Floer homology groups for $Y$ equipped with a torsion $\SpinC$
structure $\spinct$ defined in Section~\ref{HolDiskFour:sec:AbsGrade}
of~\cite{HolDiskFour}. This absolute grading has the property that if
$W$ is a cobordism from $Y_1$ to $Y_2$, (endowed with a $\SpinC$
structure $\spinc$ whose restrictions $\spinct_1$ and $\spinct_2$
respectively are both torsion), then
\begin{equation}
\label{eq:DimensionFormula}
\gr(F_{W,\spinc}(\xi))-\gr(\xi)=\frac{c_1(\spinc)^2-2\chi(W)-3\sgn(W)}{4}
\end{equation}
(c.f. Theorem~\ref{HolDiskFour:thm:AbsGrade} of~\cite{HolDiskFour}).

\begin{lemma}
\label{lemma:K3One}
For the cobordism $E-B^4$ from $S^3$ to $M\{-1,-1,-1\}$, endowed with
the $\SpinC$ structure obtained by restricting $\spinccanf$, the
generator in $\HFm_{-2}(S^3)$ is mapped to the generator
$\HFm_0(\Sigma(2,3,5))$.
\end{lemma}

\begin{proof}
For a negative-definite cobordism between integral homology spheres,
the map induced on $\HFinf$ is always an isomorphism (see
Proposition~\ref{AbsGraded:prop:NegSurgery} of~\cite{AbsGraded}).
From the dimension formula (Equation~\eqref{eq:DimensionFormula}), it
follows that the degree is raised by two.
\end{proof}

\begin{lemma}
\label{lemma:K3Two}
For the cobordism $W$
$$M\{-1,-1,-1\}\Rightarrow M\{-1,-1,0\} \Rightarrow
M\{-1,-1,1\}$$ (endowed with the $\SpinC$ structure obtained by
restricting $\spinccanf$), the induced map
$$
\begin{CD}
\Z\cong \HFm_0(M\{-1,-1,-1\})@>{F^-_{W,\spinc}}>> \HFm_{-1}(M\{-1,-1,1\})\cong \Z
\end{CD}
$$
is an isomorphism. Moreover, the if we equip $W$ with any other $\SpinC$ structure,
the induced map is trivial
\end{lemma}

\begin{proof}
Let $F$ denote the map given by the cobordism, and summing over 
all $\SpinC$ structures in $\delta H^1(M\{-1,-1,0\})$.
Dualizing (i.e. applying Theorem~\ref{HolDiskFour:thm:Duality} 
of~\cite{HolDiskFour}, and
the graded version of the duality isomorphism, c.f.
Proposition~\ref{HolDiskFour:prop:GradedDuality} of~\cite{HolDiskFour}), 
we get the following diagram:
$$
\begin{CD}
{\mathrm HF}^{-1}_-(M\{-1,1,1\})
@>{F_-}>>
{\mathrm HF}^{0}_-(M\{1,1,1\}) \\
@V{\cong}VV	@V{\cong}VV \\
\HFp_{-1}(M\{-1,-1,1\}) @>{F^+}>> \HFp_{-2}(M\{-1,-1,-1\}) \\
\end{CD}.
$$ 
From the calculations of the Floer homology groups restated above, we see
that 
$$\HFp_{-1}(M\{-1,-1,1\})\cong \Z \cong
\HFp_{-2}(M\{-1,-1,-1\}).$$ Indeed, an isomorphism is given by
composing the maps in the surgery exact sequence (see the proof of
Proposition~\ref{AbsGraded:prop:FigureEightZero} of~\cite{AbsGraded}). But this composition is precisely
$F^+$. It follows that the map $F_-$ (the map on cohomology) above is an isomorphism, and
hence (since there is no torsion present), its dual, the map 
$$F^-\colon \HFm_0(M\{-1,-1,-1\})\longrightarrow \HFm_{-1}(M\{-1,1,1\})$$
induces an isomorphism (between two groups which are isomorphic to $\Z$). 

The cobordism $W$ has $b_2(W)=2$. Indeed, we can find an embedded
torus $T_1\subset W$ which generates the image of
$H_2(M\{-1,-1,0\};\Z)$ inside $W$, and another embedded torus
$T_2\subset W$ with square zero with $T_1\cm T_2=1$. Applying the
adjunction inequality for the cobordism invariant
(Theorem~\ref{HolDiskFour:intro:Adjunction} of~\cite{HolDiskFour}) to
the embedded surface $T_2$, it follows that the only $\SpinC$
structure $\spinc\in\spinccanf+\delta H^1(M\{-1,-1,0\})$ whose associated map 
$\Fp{W,\spinc}$ is non-trivial is the restriction of $\spinccanf$ itself.
\end{proof}

\vskip.2cm
\noindent{\bf{Proof of Proposition~\ref{prop:K3Surface}}.}
According to Lemmas~\ref{lemma:K3One} and \ref{lemma:K3Two}, generator
of $\HFm_{-2}(S^3)$ is mapped to the generator of $\HFm_{-1}(M\{-1,-1,1\})\cong \Z$. 
Now, $\delta^{-1}$ of that generator is the generator of 
$\HF^+_{\red,0}(M\{-1,-1,1\})\cong \Z$. Investigating the
four exact sequences
connecting 
$$
\begin{array}{lll}
\big(M\{-1,-1,1\}, & M\{-1,0,1\}, & M\{-1,\infty,1\}\cong S^3 \big), \\
\big(M\{-1,\infty,1\}\cong S^3, & M\{-1,0,1\}, & M\{-1,1,1\}\big) ,\\
\big(M\{-1,1,1\}, & M\{0,1,1\}, & M\{\infty,1,1\}\cong S^3\big),\\
\big(M\{\infty,1,1\}\cong S^3, & M\{0,1,1\}, & M\{1,1,1\}\big),
\end{array}
$$
we see that the map
$$ \Z\cong
\HFp_{\red,0}(M\{-1,-1,1\}) \longrightarrow \HFp_{-2}(M\{1,1,1\})\cong \Z $$
induced by summing the maps induced by all $\SpinC$ structures on
the composite cobordism from $X_2-N$ is an isomorphism.
In fact, by finding square zero tori which intersect
the homology classes coming from $H_2(M\{-1,0,1\};\Z)$ and $H_2(M\{0,1,1\};\Z)$
in $X_2-N$ and applying the adjunction inequality (as in the proof
of Lemma~\ref{lemma:K3Two}), 
we see that the only $\SpinC$ structure which contributes to this sum is
the one with trivial first Chern class.
Finally, the map $$\HFp_{-2}(M\{1,1,1\})\longrightarrow
\HFp_{0}(S^3)$$ is an isomorphism (for the given $\SpinC$ structure) 
once again, 
in view of the dimension formula and the fact that $N-B^4$ 
has negative-definite intersection
form (Proposition~\ref{AbsGraded:prop:NegSurgery} of~\cite{AbsGraded}).
\qed 

%% file: nonvanish.tex
\section{The non-vanishing theorem for symplectic four-manifolds}
\label{sec:NonVanishing}

The aim of the present section is to prove
Theorem~\ref{intro:NonVanishing}. Via Donaldson's construction of
Lefschetz pencils, we will reduce this theorem to the following more
manifestly topological variant:

\begin{theorem}
\label{thm:NonVanishingLefschetz}
Let $\pi\colon X \longrightarrow S^2$ be a relatively minimal
Lefschetz fibration over the sphere with $b_2^+(X)>1$ whose generic
fiber $F$ has genus $g>1$.  Then, for the canonical
$\SpinC$ structure, we have that
\begin{eqnarray*}
\langle c_1(\spinccanf),[F]\rangle &=& 2-2g \\
\Phi_{X,\spinccanf} &=&  \pm 1.
\end{eqnarray*}
Moreover, for any other $\SpinC$ structure $\spinc\neq \spinccanf$ with
$\Phi_{X,\spinc}\neq 0$, we have that 
$$\langle c_1(\spinccanf), [F] \rangle =2-2g
< \langle c_1(\spinc),[F]\rangle.$$
\end{theorem}

One ingredient in the above proof is a related result for
three-manifolds which fiber over the circle. To state it, recall that
a three-manifold $Y$ which admits a fibration $$ \pi\colon Y
\longrightarrow S^1 $$ has a canonical $\SpinC$ structure which is
obtained as the (integrable) two-plane field, which is the kernel of
the differential of $\pi$. If $F$ is a fiber of $\pi$, then
the evaluation
$$\langle c_1(\spinccan), [F]\rangle = 2-2g.$$

\begin{theorem}
\label{thm:ThreeManifoldsFiber}
Let $Y$ be a three-manifold which fibers over the circle, with fiber genus
$g>1$, and
let $\spinct$ be a $\SpinC$ structure over $Y$ with
$$\langle c_1(\spinct),[F]\rangle = 2-2g.$$
Then, for $\spinct\neq \spinccan$, we have  that
$$\HFp(Y,\spinct)=0;$$ while
$$\HFp(Y,\spinccan)\cong \Z.$$
\end{theorem}

Indeed, we also establish the following result, which 
bridges the above two theorems:

\begin{theorem}
\label{thm:LefschetzDisk}
Let $\pi\colon W \longrightarrow \CDisk$ be a relatively minimal
Lefschetz fibration over the disk with fiber genus $g>1$,
and let $Y=-\partial W$.
Then, there is a unique $\SpinC$ structure 
$\spinc$ over $W$ for which
$$\langle c_1(\spinc),[F]\rangle = 2-2g,$$
and the induced map 
$$F^+_{W-B^4,\spinc}\colon \HFp(Y,\spinc|Y)\longrightarrow \HFp(S^3)$$
is non-trivial; and that is the canonical $\SpinC$ structure $\spinccanf$.
Indeed, the induced map
$$F^+_{W-B^4,\spinccanf}\colon \HFp(Y,\spinccanf|Y)\longrightarrow \HFp_0(S^3)\cong
\Z$$
is an isomorphism.
\end{theorem}

We prove the above three theorems, in reverse order.

In fact, we prove several special cases of these theorems first.  It
will be convenient to fix some notation. Suppose that $W$ is
some four-manifold which admits a Lefschetz fibration $\pi$ (over some
two-manifold possibly with boundary). Then we let $${\mathfrak
S}(W)=\{\spinc\in\SpinC(W)\big| \langle c_1(\spinct),[F]\rangle
=2-2g\}.$$ (This is a slight abuse of notation: ${\mathfrak S}(W)$
depends on the Lefschetz fibration $\pi$, not just the four-manifold
$W$.) Similarly, if $Y$ is a three-manifold which fibers over the
circle, we let
$${\mathfrak T}(Y)=\{\spinct\in\SpinC(Y)\big| \langle
c_1(\spinct),[F]\rangle =2-2g\}.$$
We will also let $\HFp(Y,{\mathfrak T}(Y))$ denote the direct sum
$$\HFp(Y,{\mathfrak T}(Y))=\bigoplus_{\spinct\in{\mathfrak T}(Y)}\HFp(Y,\spinct).$$

\begin{lemma}
\label{lemma:MaxSpinCIndepOfFib}
Let $\pi\colon W\longrightarrow [1,2]\times S^1$ be a relatively minimal 
Lefschetz
fibration with fiber genus $g>1$ 
over the annulus, which connects a pair of three-manifolds
$Y_1$ and $Y_2$ (which fiber over the circle), then for some choice of signs,
the map
$$\sum_{\spinc\in{\mathfrak S}(W)} \pm \Fp{W,\spinc}\colon
\HFp(Y_1,{\mathfrak T}(Y_1))
\longrightarrow
\HFp(Y_2,{\mathfrak T}(Y_2))
$$
induces an isomorphism.
\end{lemma}

\begin{proof}
Note that whereas ${\mathfrak S}(W)$ can easily by infinite;
according to the finiteness properties  for 
the maps associated to cobordisms (Theorem~\ref{HolDiskFour:thm:Finiteness}
of~\cite{HolDiskFour})
there are only finitely many $\spinc\in{\mathfrak S}(W)$ for
which $\Fp{W,\spinc}$ is non-trivial. 

First assume that the Lefschetz fibration $\pi$ has a single node.  In
this case, $W$ can be viewed as the cobordism obtained by attaching a
single two-handle to $Y=Y_1$ along a curve $K$ in the fiber of $\pi$,
with framing $-1$ (with respect to framing $K$ inherits from the fiber
$F\subset Y$); in particular, $Y_2=Y_{-1}(K)$. Moreover, since the
Lefschetz fibration is relatively minimal, the curve $K$ is
homotopically non-trivial as a curve in $F$.  Now, if $Y_0(K)$ is the
three-manifold obtained as zero-surgery along $K$ then the cobordism
from $Y$ to $Y_0$ also maps to the circle (by a map $\pi_0$ which is
no longer a fibration, but which extends the map $\pi$ from $Y$ to
$S^1$).  Clearly, if $\spinc$ is any $\SpinC$ structure which extends
over $W_0$, the restriction of $c_1(\spinc)$ to a generic fiber of
$\pi_0\colon Y_0(K)\longrightarrow S^1$ is also $2-2g$. However, since
$K$ is homotopically non-trivial, the Thurston norm of the homology
class of this fiber in $Y_0(K)$ is smaller than $2-2g$, so the
adjunction inequality for $\HFp$
(Theorem~\ref{HolDiskTwo:thm:Adjunction} of~\cite{HolDiskTwo}) ensures
that $\HFp(Y_0,\spinc|Y_0)=0$. Thus, the lemma follows immediately
from the surgery long exact sequence for $\HFp$ (see
Theorem~\ref{HolDiskTwo:thm:GeneralSurgery} of~\cite{HolDiskTwo}):
$$ ...\longrightarrow
\HFp(Y,{\mathfrak T}(Y)) \longrightarrow \HFp(Y_{-1}(K),{\mathfrak
T}(Y_{-1}(K)))
\longrightarrow \HFp(Y_0(K),{\mathfrak T}(Y_0))=0 \longrightarrow...
$$ In the above sequence,  ${\mathfrak T}(Y_0)$ denotes those
$\SpinC$ structures whose evaluation on the homology class of a fiber of $\pi_0$
(which is no longer a fibration) is given by
$2-2g$, where now $g$ still denotes  the genus of the fibration for $Y$. 

The case of multiple nodes follows immediately by the composition law.
\end{proof}

\begin{lemma}
\label{lemma:MaxSpinC}
If $\pi\colon Y \longrightarrow S^1$ is a surface bundle over $S^1$,
with fiber genus $g>1$, then there is a unique $\SpinC$
structure $\spinct\in{\mathfrak T}(Y)$ with 
$\HFp(Y,\spinct)\neq 0$. In fact, 
$$\HFp(Y,\spinct)\cong \Z.$$
\end{lemma}

\begin{proof}
Note that the mapping class group is generated as a monoid by (right-handed)
Dehn twists. This is equivalent to the claim that if 
$p_1\colon Y_1\longrightarrow S^1$ and
$p_2\colon Y_2\longrightarrow S^1$ 
any two fibrations over the circle whose fiber has the same genus, then
we can extend the two fibrations to form a 
relatively Lefschetz fibration over the annulus. It follows from 
Lemma~\ref{lemma:MaxSpinCIndepOfFib} that for a genus $g$ fibration over the circle
$\HFp(Y,{\mathfrak T}(Y))$ is independent of the monodromy map, and depends only
on the genus $g$. 

Thus, for each $g>1$, it suffices to find some fibered three-manifold
for which the lemma is known to be true. 
For this purpose, let $Y=Y(g)$ be the
zero-surgery on the torus knot $K$ of type $(2,2g+1)$. This is a
fibered three-manifold whose fiber has genus $g$. 
Writing the symmetrized Alexander polynomial of $K$
as
$$\Delta_K(T)=-\sum_{i=-g}^g(-T)^i=a_0+\sum_{i=1}^d a_i(T^{i}+T^{-i})$$ 
it is shown in Proposition~\ref{AbsGraded:prop:CorrTermTorusKnots}
of~\cite{AbsGraded} that if $\spinct$ is a $\SpinC$ structure
over $Y$ with 
$$\langle
c_1(\spinct), [F]\rangle = 2i\neq 0,$$ then
$\HFp(Y,\spinct)$ is a free Abelian group of rank
$$\sum_{j=1}^\infty j a_{|i|+j}.$$
In particular, when $\langle c_1(\spinct),[F]\rangle =2-2g$, it follows immediately that
$\HFp(Y,\spinct)\cong \Z$.
\end{proof}

\begin{lemma}
\label{lemma:SigmaTimesDisk}
Let $\Fiber$ be an oriented surface of genus $g>0$, and consider the
cobordism $W$ from $S^3$ to $\Fiber\times S^1$ obtained by puncturing
the product $\Fiber\times D^2$ in a single point.  Let $\spinccanf$
denote the $\SpinC$ structure over $W$ with $\langle
c_1(\spinccanf),[\Fiber]\rangle = 2-2g$. Then, the induced map
$$\Fp{W,\spinccanf}\colon
\HFp(\Fiber\times S^1,\spinccan)\longrightarrow
\HFp_0(S^3)\cong \Z$$
is an isomorphism, as is the induced map
$$(1-U^{g-1})\Fm{W,\spinccanf}\colon \Z\cong \HFm_{-2}(S^3)\longrightarrow
\HFmRed(\Fiber\times S^1,\spinccan)\cong \Z
\subset \HFm(\Fiber\times S^1,\spinccan).$$ 
\end{lemma}

\begin{proof}
To see the claim about $\Fp{W,\spinccanf}$, it suffices to embed the
cobordism $(W,\spinccanf)$ into a closed four-manifold $(X,\spinc)$ 
with $b_2^+(X)>1$, so
that $\spinc|W=\spinccanf$ and $\Phi_{X,\spinc}=\pm 1$.  To see why
this suffices, observe that $U\cdot \HFp(F\times S^1,\spinccan)=0$, so
$\Fp{W,\spinc}$ must take $\HFp(F\times S^1,\spinccan)$ into
$\HFp_0(S^3)\cong \Z$. In general, the image of such a map 
consists of multiples of some integer $d$. Now,
take an admissible cut of $X=X_1\#_N X_2$ which is disjoint from $F$,
and so that $F\subset X_2$ (such a cut is found by taking any embedded
surface $\Sigma$ of positive square which is disjoint from $F$). It
then follows that for each $\SpinC$ structure $\spinc\in\SpinC(X)$
which restricts to $W$ as $\spinccanf$, the sum of invariants
$$\sum_{n\in\Z}\Phi_{X,\spinc+\PD[\Sigma]}$$ is divisible by $d$.
In fact, it is a straightforward consequence of the
dimension formula that the part of this sum which is homogeneous of degree zero is
the invariant $\Phi_{X,\spinc}$, and this, in turn, forces $d=\pm 1$,
so that the claimed map is an isomorphism.

Now, such four-manifolds can be found for all possible genera $g$ in
the blow-ups of the $K3$ surface, in light of the blow-up formula and
the $K3$ calculation. Specifically, for each genus $g$, we can find an
embedded surface $\Sigma\subset K3$ with $\Sigma\cm
\Sigma=2g-2$, 
for instance, by taking a single section of an elliptic fibration of
$K3$, which is a sphere of self-intersection number $-2$, and
attaching $g$ copies of the fiber. In the $2g-2$-fold blow-up,
$\Sigma$ has a proper transform ${\widehat \Sigma}$ with
${\widehat\Sigma}\cm{\widehat \Sigma}=0$. Consider the $\SpinC$
structure ${\widehat\spinc}$ with
$c_1({\widehat\spinc})=-\PD[E_1]-...-\PD[E_{2g-2}]$, so that $\langle
c_1({\widehat\spinc}),[{\widehat\Sigma}]\rangle=2-2g$; i.e. the
tubular neighborhood of $\Sigma$ is $W$, and ${\widehat \spinc}$ is
an extension of $\spinccanf$. According to Proposition~\ref{prop:K3Surface}
and the blowup formula (Theorem~\ref{HolDiskFour:intro:BlowUp}
of~\cite{HolDiskFour}), $\Phi_{{\widehat X},{\widehat\spinc}}=\pm 1$.

The statement about $\HFm$ follows similarly, by choosing the cut for
$X$ so that the surface $F$ lies in $X_1$.
\end{proof}

\begin{lemma}
\label{lemma:NonSeparatingCase}
Let $\pi\colon W \longrightarrow \CDisk$ be a Lefschetz fibration over
the disk, whose singular fibers are all non-separating nodes. Then,
$\pi\colon W \longrightarrow \CDisk$ can be embedded in a Lefschetz
fibration $V$ over a larger disk with the property that the canonical $\SpinC$ structure
$\spinccanf$ is the only $\SpinC$ structure in $\spinc\in {\mathfrak S}(V)$ 
for which 
$$\Fp{V,\spinc}\colon \HFp(\partial V,{\mathfrak T}(\partial V)) 
\longrightarrow
\HFp_0(S^3)\cong \Z$$
is non-trivial; and indeed, $\Fp{V,\spinccanf}$ is an isomorphism.
\end{lemma}

\begin{proof}
We claim that any Lefschetz fibration over the disk with
non-separating fibers can be embedded into a Lefschetz fibration over
the disk with nodes corresponding to (isotopic translates) of
the standard curves $\{\tau_1,...,\tau_m\}$
described in Theorem~\ref{thm:GenMappingClassGroup}.
This is constructed as follows. Suppose that $W$ is described by monodromies
which are Dehn twists around curves $(C_1,...,C_n)$. Then, we can find 
automorphisms of $F$, $\phi_1$,...,$\phi_n$, so that $\phi_i(\tau_1)=C_i$.
We then express each $\phi_i=D(\tau_{m_{i,1}})\cm...
\cm D(\tau_{m_{i,\ell_i}})$.
We let $V$ be the Lefschetz fibration over the disk with monodromies obtained
by juxtaposing $\tau_{m_{i,1}},...,\tau_{m_{i,\ell_i}},\tau_1$ for $i=1,...,n$,
union as many $\tau_i$ as it takes to span all of $H_1(\Sigma;\Z)$.
By performing Hurwitz moves, we obtain a subfibration with monodromies
$(\phi_1(\tau_1),...,\phi_n(\tau_1))$; i.e. we have embedded $W$ in $V$.

Next, we argue that $V$ has the required form.  According to
Lemmas~\ref{lemma:MaxSpinC}, \ref{lemma:SigmaTimesDisk}, and
\ref{lemma:MaxSpinCIndepOfFib}, we see that
$$\sum_{\spinc\in{\mathfrak S}(V)} \Fp{V-B^4,\spinc}\colon
\HFp(\partial V,\spinct)\cong \Z \longrightarrow \HFp_{-2}(S^3)$$
is an isomorphism. We claim that $\spinccanf$ is the only $\SpinC$ structure
in the sum with non-zero contribution.

Note that $H_1(V;\Z)$ is the quotient of $\Z^{2g}$ by the homology
homology classes of the vanishing cycles for $V$, so we have arranged
that $H_1(V;\Z)=0$; in particular, $H^2(V;\Z)$ has no torsion.  It
follows that the $\SpinC$ structure $\spinccanf$ is uniquely
determined by the evaluation of its first Chern class on the various
two-dimensional homology classes in $V$.  Moreover, if we choose the
translates of the various $\tau_i$ carefully, so that parallel copies
of the same $\tau_i$ remain disjoint, then we can find a basis for
$H_2(V;\Z)$ consisting of $[F]$ and surfaces ${\widehat P}$ obtained
by ``capping off'' submanifolds-with-boundary $P\subset F$ whose
boundaries consist of copies of the vanishing cycles. Suppose, next,
that ${\widehat P}_1$ is induced from a relation $P_1$ in $F$ with
this form, and let $m$ denote the number of its boundary
components. Then, the relation $F-P_1=P_2$ also has this form (and has
the same number of boundary components), and its closed extension
${\widehat P}_2$ satisfies the following elementary properties (see
Lemma~\ref{lemma:PeriodicDomainsInFibration}):
\begin{eqnarray*}
[F]&=&[{\widehat P}_1] + [{\widehat P}_2], \\
g(F)&=&g({\widehat P}_1)+g({\widehat P}_2) +m-1, \\
m&=& -[{\widehat P}_1]^2=-[{\widehat P}_1]^2 = 
[{\widehat P}_1] \cm[{\widehat P}_2]
\end{eqnarray*}

Now suppose that $\spinc\in{\mathfrak S}(V)$ is a $\SpinC$ structure 
for which $\Fp{W,\spinc}$ is non-trivial. Then, the above equations, and the condition that $\langle c_1(\spinc),[F]\rangle = 2-2g$ say that
$$\left(\langle c_1(\spinc), [{\widehat P}_1]\rangle - [{\widehat P}_1]\cm [{\widehat P}_1]\right) + 
\left(\langle c_1(\spinc), [{\widehat P}_2]\rangle - [{\widehat P}_2]\cm [{\widehat P}_2]\right) = \left(2-2g({\widehat P}_1)\right) + 
\left(2-2g({\widehat P}_2)\right).
$$
Now, either 
$$\langle c_1(\spinc), [{\widehat P}_1]\rangle - [{\widehat P}_1]\cm [{\widehat P}_1] = 2-2g([{\widehat P}_1]),$$
in which case (according to Lemma~\ref{lemma:PeriodicDomainsInFibration}),
\begin{equation}
\label{eq:EvaluateCOneSpinC}
\langle c_1(\spinc), [{\widehat P}_1] \rangle = 
\langle c_1(\spinccanf),[{\widehat P}_1] \rangle,
\end{equation}
or, after possibly switching the roles of ${\widehat P}_1$
and ${\widehat P}_2$,  we have that
\begin{equation}
\label{ineq:AdjunctionInequality}
\langle c_1(\spinc),[{\widehat P}_1]\rangle -[{\widehat P}_1]
\cm[{\widehat P}_1]
\leq -2g([{\widehat P}_1]).
\end{equation}

Inequality~\eqref{ineq:AdjunctionInequality} 
is ruled out by the adjunction relation,
Theorem~\ref{thm:AdjunctionRelation}, as follows. By adding trivial
two-handles to ${\widehat P}_1$ if necessary, we obtain an
embedded surface with $\langle c_1(\spinc),[\Sigma]\rangle =
-2g+\Sigma\cm\Sigma$. There are two cases, according to whether
$g(\Sigma)=0$ or $g(\Sigma)>0$. In the latter case,
the adjunction relation gives some word $\xi(\Sigma)$ of degree
$2g(\Sigma)>0$ in $\Alg(\Sigma)$ with the property that $$
\Fp{V,\spinc}(\cdot)=\Fp{V,\spinc+\epsilon\PD[\Sigma]}(\xi(\Sigma)
\otimes \cdot).$$
Observe that homology classes in $\Sigma$ 
are all homologous to classes in the fiber $F$ in $\partial V$, so the 
action by $\xi(\Sigma)$ appearing above can be interpreted as the
action by an element of positive degree in $\Z[U]\otimes_{\Z}\Wedge^*(H_1(Y)/\Tors)$ on $\HFp(\partial V,\spinccan)$.
But all such elements annihillate $\HFp(\partial V,\spinccan)$ (since it
is supported in a single dimension). Thus,
the only remaining possibility is that $g(\Sigma)=0$, in which case no
handles were added to ${\widehat P}_1$. In thise case, the adjunction
relation ensures that the $\SpinC$ structure $\spinc-\PD[{\widehat
P}_1]$ has non-trivial invariant, while $$\langle
c_1(\spinc),[{\widehat P}_2]\rangle - [{\widehat P}_2]\cm [{\widehat
P}_2]=4-2g({\widehat P}_2).$$ But then, $$\langle
c_1(\spinc-\PD[{\widehat P}_1]),[{\widehat P}_2]\rangle - [{\widehat
P}_2]\cm [{\widehat P}_2]=4-2g({\widehat P}_2)-2m.$$ Next, observe
that $m>1$, since the vanishing cycles for $V$ are all homotopically
non-trivial. Moreover, if $m=2$, then $g({\widehat P}_2)=g(F)$.
Thus, using ${\widehat P}_2$ in place of ${\widehat P}_1$, and 
$\spinc-\PD[{\widehat P}_1]$ in place of $\spinc$,
we obtain the same contradiction as before.

The contradiction to Inequality~\eqref{ineq:AdjunctionInequality}
leads to the conclusion that Equation~\eqref{eq:EvaluateCOneSpinC}
holds for all choices of ${\widehat P}_1$.  But these surfaces,
together with $[F]$, generate the homology of $V$. Thus, we have shown
that $\spinc=\spinccanf$, as claimed.
\end{proof}

\vskip.2cm
\noindent{\bf{Proof of Theorem~\ref{thm:ThreeManifoldsFiber}.}}
According to Lemma~\ref{lemma:MaxSpinC}, there is a unique
$\spinct\in{\mathfrak T}(Y)$ with $\HFp(Y,\spinct)\neq 0$, and for
$\spinct$, we have that $\HFp(Y,\spinct)\cong \Z$. It remains to
identify $\spinct$ with the canonical $\SpinC$ structure.  As in the
proof of the lemma, we constructed a Lefschetz fibration over the
annulus which connects $Y$ with $S^1\times \Sigma$. By attaching
$\CDisk\times\Sigma$ to the $S^1\times \Sigma$ boundary component, we
obtain a Lefschetz fibration $W$ over the disk. Indeed, since the
mapping class group is generated by Dehn twists along
non-separating curves, we can choose $W$ so that
Lemma~\ref{lemma:NonSeparatingCase} applies to $W$. In particular, in
this case, the canonical $\SpinC$ structure $\spinccanf$ in ${\mathfrak S}(W)$ 
induces a non-trivial map
$\Fp{W,\spinc}$. The result follows, since $\spinccanf|Y=\spinccan$.
\qed

\begin{lemma}
\label{lemma:SeparatingCase}
Let $W$ be a relatively minimal Lefschetz fibration over the annulus, all of whose nodes
are separating. Then, the only $\SpinC$ structure $\spinc\in{\mathfrak
S}(W)$ for which the map
$$\Fp{W,\spinc}\colon \HFp(Y_1,\spinc|Y_1)\cong \Z
\longrightarrow \HFp(Y_2,\spinc|Y_2)\cong \Z$$
is non-trivial is the canonical $\SpinC$ structure. And for that $\SpinC$ structure,
the induced map is an isomorphism.
\end{lemma}

\begin{proof}
According to Lemmas~\ref{lemma:MaxSpinC}, \ref{lemma:SigmaTimesDisk},
and \ref{lemma:MaxSpinCIndepOfFib}, we see that
$$\sum_{\spinc\in{\mathfrak S}(W)} \Fp{W,\spinc}\colon
\sum_{\spinct\in{\mathfrak T}(Y)}
\HFp(Y,\spinct) \longrightarrow 
\sum_{\spinct\in{\mathfrak T}(S^1\times F)}
\HFp(S^1\times F,\spinct)\cong \HFp(S^1\times F,\spinccan)\cong \Z$$
is an isomorphism. 

Now, observe that $W$ is a cobordism which is obtained by attaching a
sequence of two-handles along null-homologous curves. Thus, a $\SpinC$
structure over $W$ is uniquely characterized by its restriction to one
of its boundary components, and its evaluations on the two-dimensional
homology classes introduced by the two-handles. According to
Theorem~\ref{thm:ThreeManifoldsFiber}, the restriction to the boundary
must agree with the canonical $\SpinC$ structure. Each node has, as
fiber, a union of two surfaces meeting at a point: i.e. we obtain a
pair of embedded surfaces $g({\widehat P}_1)+g({\widehat P}_2)=g(F)$
and ${\widehat P}_1^2={\widehat P}_2^2=-1$. Moreover, since the
fibration is assumed to be relatively minimal, $g({\widehat P}_1)>0$
and $g({\widehat P}_2)>0$. Thus, applying the adjunction relation as
in the proof of Lemma~\ref{lemma:NonSeparatingCase}, we see that
$$\langle c_1(\spinc),[{\widehat P}_1] \rangle = \langle
c_1(\spinccanf),[{\widehat P}_1] \rangle.$$
It is easy to see that the homology classes of the form
$[{\widehat P}_1]$ (one for each node) generate
$H_2(W;\Z)/H_2(Y;\Z)$. Thus, it follows that $\spinc=\spinccanf$.
\end{proof}

\noindent{\bf{Proof of Theorem~\ref{thm:LefschetzDisk}.}}
Let $\pi\colon X\longrightarrow \CDisk$ be the Lefschetz fibration.
By combining Lemma~\ref{lemma:MaxSpinCIndepOfFib}, 
Lemma~\ref{lemma:MaxSpinC}, and Lemma~\ref{lemma:SigmaTimesDisk}, we see that the map
$$\sum_{\spinc\in{\mathfrak S}(X)} \Fp{X,\spinc}\colon
\bigoplus_{\spinct\in{\mathfrak T}(Y)}\HFp(Y,\spinct)
\longrightarrow
\HFp_0(S^3)$$
induces an isomorphism. 

We can find a subdisk $\CDisk_0\subset \CDisk$ which contains all the
fibers with non-separating nodes. Let $X_0\subset X$ denote its
preimage. According to Lemmas~\ref{lemma:MaxSpinCIndepOfFib},
\ref{lemma:MaxSpinC}, and \ref{lemma:SigmaTimesDisk}, there must
be at least one $\SpinC$ structure $\spinc\in{\mathfrak S}(X)$ for
which the map $$\Fp{X,\spinc}\colon \HFp(Y,\spinc|Y)\cong \Z
\longrightarrow
\HFp(S^3)$$
is non-trivial. According to Lemma~\ref{lemma:NonSeparatingCase}, its
restriction $\spinc|X_0$ is the canonical $\SpinC$ structure; according
to Lemma~\ref{lemma:SeparatingCase}, its restriction $\spinc|X-X_0$ is
also the canonical $\SpinC$ structure. Now, the map
$H^1(X-X_0)\longrightarrow H^1(Y;\Z)$ is an isomorphism, since $X-X_0$
is obtained from $Y\times[0,1]$ by attaching two-handles along
null-homologous curves. Thus, the only $\SpinC$ structure whose
restrictions to both $X-X_0$ and $X_0$ agree with $\spinccanf$ is the
canonical $\SpinC$ structure $\spinccanf$ itself.
\qed
\vskip.2cm

\vskip.2cm
\noindent{\bf{Proof of Theorem~\ref{thm:NonVanishingLefschetz}.}}
We decompose $X=X_1\#_{S^1\times \Sigma_g} X_2$
where $X_1$ is the pre-image of a disk in the Lefschetz fibration which
contains no singular points (in particular, $X_1=\CDisk\times F$).
According to Proposition~\ref{prop:OtherCut}, 
$$\Fp{W_2,\spinc_2}\circ
\PiRed_N\circ
\Fm{W_1,\spinc_1}=
\sum_{\{\spinc\in\SpinC(W)|\spinc|W_1=\spinc_1,\spinc|W_2=\spinc_2\}}
\Fmix{W,\spinc}.$$ 
Now, by Lemma~\ref{lemma:SigmaTimesDisk}, 
$$\PiRed\circ\Fm{W_1,\spinc_1}\colon \Z\cong \HFm_{-2}(S^3)\longrightarrow
\HFpRed(S^1\times \Sigma_g)\cong \Z$$
is an isomorphism. Similarly, according to Theorem~\ref{thm:LefschetzDisk},
$$\Fp{W_2,\spinc_2}\colon 
\HFpRed(S^1\times \Sigma_g)\cong \Z\longrightarrow \HFp_0(S^3)\cong \Z$$
is an isomorphism. Thus, we conclude that $$1 = \sum_{\eta\in \delta
H^1(\Sigma\times S^1)} \pm\Phi_{X,\spinc+\eta}.$$ Observe, however,
that $\delta H^1(\Sigma\times S^1)$ is one-dimensional; 
in fact, 
the $\SpinC$ structures in the $\delta H^1(\Sigma\times S^1)$-orbit
are of the form $\spinccanf+\Z \PD[F]$. By the dimension formula, the only such
$\SpinC$ structure which has degree zero is $\spinccanf$ (using the
adjunction formula and the fact that the fiber genus $g>1$). If $k>0$
we see that $\Fmix{W,\spinc-k\PD[F]}$ is zero. If
$\Fmix{W,\spinc+k\PD[F]}$ were non-zero, the expression
$\Fp{W_2,\spinc_2}\circ
\PiRed_N\circ
\Fm{W_1,\spinc_1}(U)$ would have to be non-zero.
But this is impossible, since
$U$ annihilates $\HFpRed(S^1\times \Sigma_g,\spinccan)$.  

Finally, we observe that the usual adjunction inequality for
surfaces with square zero (Theorem~\ref{HolDiskFour:intro:Adjunction} of~\cite{HolDiskFour})
ensures that if
$$\langle c_1(\spinccanf), [F] \rangle =2-2g
> \langle c_1(\spinc),[F]\rangle,$$
then $\Phi_{X,\spinc}\equiv 0$.
\qed
\vskip.2cm

\vskip.2cm
\noindent{\bf{Proof of Theorem~\ref{intro:NonVanishing}.}}
First, observe that the conditions on $\omega$ in Theorem~\ref{intro:NonVanishing}
are all open conditions, so it suffices to prove the theorem in the case
where $\omega$ has rational periods.
According to Donaldson's theorem,
any sufficiently large multiple $N\omega$ gives rise to a Lefschetz
pencil. Specifically, if we blow up $X$ sufficiently many times, we
get a new symplectic manifold $({\widehat X},{\widehat \omega})$ with
the property that $$N\omega-\sum_{i=1}^m \PD[E_i]$$ is
Poincar\'e dual to the fiber of a Lefschetz fibration over $S^2$. Here,
$\{E_i\}_{i=1}^m$ are the exceptional spheres in ${\widehat X}$.
In particular, for any $\SpinC$ structure $\spinc\in\SpinC(X)$, we have that
\begin{equation}
\label{eq:RelateWithVolume}
\langle c_1({\widehat\spinc}),[F]\rangle = 
\langle c_1(\spinc),N\omega\rangle - m.
\end{equation}

Clearly, the canonical $\SpinC$ structure of $({\widehat
X},{\widehat\omega})$ is the blow-up of the canonical $\SpinC$
structure of $(X,\omega)$, so according to the blow-up formula for
$\Phi$, follows that $\Phi_{X,\spinccan}=\pm 1$ if and only if
$\Phi_{{\widehat X},{\widehat \spinccan}}=\pm 1$. But the latter
equation follows, according to
Theorem~\ref{thm:NonVanishingLefschetz}.  

For suitable choice of $N$, we can arrange for the Lefschetz fibration
to be relatively minimal, see~\cite{Smith} and~\cite{AurouxKatz}.  In
this case, if $\spinc\in\SpinC(X)$ is any structure with
$\Phi_{X,\spinc}\neq 0$, then its blowup ${\widehat\spinc}$ satisfies
$\Phi_{{\widehat X},{\widehat\spinc}}\neq 0$. Thus, the inequality
stated in this theorem is equivalent to the corresponding inequality
from Theorem~\ref{thm:NonVanishingLefschetz}, in view of
Equation~\eqref{eq:RelateWithVolume}.
\qed

%% file: sympthom.tex
\section{The genus-minimizing properties of symplectic submanifolds}
\label{sec:SymplecticThom}

In the case where $b_2^+(X)>1$, Theorem~\ref{intro:SympThom} is now an
easy consequence of Theorem~\ref{thm:AdjunctionRelation} and
Theorem~\ref{intro:NonVanishing}. For this implication, we
follow~\cite{SympThom}

\vskip.2cm
\noindent{\bf{Proof of Theorem~\ref{intro:SympThom} when $b_2^+(X)>1$.}}
If the theorem were false, we could find a symplectic manifold
$(X,\omega)$ and a pair $\Sigma, \Sigma'\subset X$ of homologous,
smoothly-embedded submanifolds, with $\Sigma$ symplectic, and
$g(\Sigma')<g(\Sigma)$. By blowing up $X$ and taking the proper
transform of $\Sigma$ as necessary, we can assume that $\langle
c_1(\spinccanf),[\Sigma]\rangle <0$.  By attaching handles to $\Sigma'$
as necessary, we can arrange for $g(\Sigma')=g(\Sigma)-1$.  Then, the
adjunction formula for $\Sigma$ gives us that $$ \langle
c_1(\spinccanf),[\Sigma'] \rangle - [\Sigma']\cm [\Sigma'] =
-2g(\Sigma'). $$ Theorem~\ref{intro:NonVanishing} says that
$\Phi_{X,\spinccanf}$ is non-trivial, so according to
Theorem~\ref{thm:AdjunctionRelation}, 
$\Phi_{X,\spinccanf-\PD[\Sigma']}$ is non-trivial, as well. But since
$$\langle \omega,c_1(\spinccanf-\PD[\Sigma'])\rangle = \langle
\omega,c_1(\spinccanf) \rangle - 2\langle \omega,[\Sigma]\rangle < 
\langle
\omega,c_1(\spinccanf) \rangle,$$
we obtain the desired contradiction to Theorem~\ref{intro:NonVanishing}.
\qed
\vskip.2cm

For the case where $b_2^+(X)=1$, we appeal directly to the analogue of
Theorem~\ref{thm:NonVanishingLefschetz}. 

Specifically, recall that if $\pi\colon X\longrightarrow S^2$ is a 
Lefschetz fibration with genus $g>1$, then
$\langle c_1(\spinccanf),[F]\rangle=2-2g\neq 0$, so we have an invariant
$\Phi_{X,\spinc,L}$ in the sense of Subsection~\ref{subsec:BTwoOne},
where $L$ is the line containing $F$ in $H^2(X;\Q)$.
The proof of Theorem~\ref{thm:NonVanishingLefschetz} gives:

\begin{theorem}
\label{thm:NonVanishingLefschetzBOne}
Let $\pi\colon X \longrightarrow S^2$ be a relatively minimal
Lefschetz fibration over the sphere with $b_2^+(X)=1$ whose generic
fiber $F$ has genus $g>1$.  Then, for the canonical
$\SpinC$ structure, we have that
\begin{eqnarray*}
\langle c_1(\spinccanf),[F]\rangle &=& 2-2g \\
\Phi_{X,\spinccanf,L} &=&  \pm 1,
\end{eqnarray*}
where $L$ denotes the line in $H^2(X;\Q)$ containing $[F]$.
Moreover, for any other $\SpinC$ structure $\spinc\neq \spinccanf$ with
$\Phi_{X,\spinc,L}\neq 0$, we have that 
\begin{equation}
\label{ineq:OtherGuysInequality}
\langle c_1(\spinccanf), [F] \rangle =2-2g
< \langle c_1(\spinc),[F]\rangle.
\end{equation}
\end{theorem}

\vskip.2cm
\noindent{\bf{Proof of Theorem~\ref{intro:SympThom} when $b_2^+(X)=1$.}}
Once again, if the theorem were false, we would be able to find
homologous surfaces $\Sigma$ and $\Sigma'$ in $(X,\omega)$ with
$\Sigma$ symplectic and $g(\Sigma')=g(\Sigma)-1$. We claim that for
sufficiently large $N$, we can find a relatively minimal Lefschetz
fibration on some blowup ${\widehat X}$ whose fiber $F$ satisfies
$F\cm {\widehat \Sigma}=0$, where ${\widehat \Sigma}$ is some suitable
proper transform of $\Sigma$. Specifically, if $\omega\cm \Sigma=c$
(which we can assume is an integer), then provided that $N\omega^2>c$, we can
let ${\widehat \Sigma}$ represent the homology class 
$$	[{\widehat\Sigma}]=[\Sigma]-[E_1]-...-[E_{Nc}]	$$ 
inside the Lefschetz fibration
obtained by blowing up the Lefschetz pencil for $N\omega$. The
homology class of the fiber here is given by
$$[F]=N[\omega]-[E_1]-...-[E_{M}],$$ where $M=N^2\omega^2$. Of course,
Theorem~\ref{thm:NonVanishingLefschetzBOne} ensures that 
$\Phi_{X,\spinccanf,L}\not\equiv 0$.
We can then find
a new embedded surface $F'$ representing $F$, but which is disjoint
from $\Sigma'$, and cut $X$ along $F'\times S^1$ into two pieces, one
of which is a tubular neighborhood of $F'$. For this cut,
Theorem~\ref{thm:AdjunctionRelation} shows that 
$\Phi_{X,\spinccanf\pm\PD[\Sigma],L}$ is also non-trivial. But
since 
$$\langle c_1(\spinccanf\pm \PD[\Sigma]),[F]\rangle =
\langle c_1(\spinccanf),[F]\rangle,$$
this violates Inequality~\eqref{ineq:OtherGuysInequality}.
\qed
\vskip.2cm

%% file: someplumbs.tex
\section{A class of three-manifolds with $\HFpRed(Y)=0$}
\label{sec:SomePlumbings}

We now prove the following:

\begin{theorem}
\label{thm:FloerHomology}
Let $Y$ be a three-manifold which can be obtained as a 
plumbing of spheres specified by a weighted graph $(G,m)$ which
satisfies the following conditions:
\begin{itemize}
\item $G$ is a disjoint union of trees
\item at each vertex in $G$, we have that
\begin{equation}
\label{eq:VertexInequality}
m(v)\geq d(v).
\end{equation}
Then, $\HFpRed(Y)=0$.
\end{itemize} 
\end{theorem}

Note that any lens space can be expressed as a plumbing of two-spheres
along a graph $(G,m)$ satisfying the above hypotheses. (Indeed, the
graph is linear: it is connected, each vertex has degree at most two,
and multiplicity at least two.)

Any Seifert fibered space $Y$ with $b_1(Y)\leq 1$ and which is not a
lens space is obtained as a plumbing along a star-like graph: the
graph is connected, has a unique vertex (the ``central node'') with
degree $n>2$, and all other vertices have degree at most two and multiplicity at least
two. The degree of the central node agrees with the number of
``singular fibers'' of the Seifert fibration, and its multiplicity $b$
is one of the Seifert invariants of the fibration. Thus, a Seifert
fibration satisfies the hypotheses of the above theorem when $b\geq
n$.

\begin{remark}
An easy inductive argument similar to the proof given below also gives
the relative grading. Suppose that $(G,m)$ is a weighted graph satisfying the hypotheses of
Theorem~\ref{thm:FloerHomology}, with the additional hypothesis that
$Y=-Y(G,m)$ is a rational homology three-sphere (this in turn is equivalent to the 
hypothesis that each component of $G$ contains at least one vertex for which
Inequality~\eqref{eq:VertexInequality} is strict), and let $W(G,m)$ be 
the four-manifold obtained by plumbing two-sphere
bundles according to a weighted graph $(G,m)$, and let $W=-W(G,m)$ be the 
plumbing with negative-definite intersection form.
Then for each $\spinct\in \SpinC(Y)$, letting ${\mathfrak K}(\spinct)$
denote the set of characteristic vectors $K\in H^2(W;\Z)$ for which $K|Y=c_1(\spinct)$,
we have that
\begin{equation}
\label{eq:DCalc}
d(Y,\spinct)=\min_{K\in{\mathfrak K}(\spinct)}
\frac{K^2+|G|}{4},
\end{equation}
where $|G|=\rk(H_2(W))$ denotes the number of vertices in $G$.
Indeed, Equation~\eqref{eq:DCalc}
remains true even in the case where the graph has a single vertex
where Inequality~\eqref{eq:VertexInequality} fails, which includes all
Seifert fibered rational homology three-spheres.
We return to these topics, and the more general issue of determining $\HFp$ for
trees with arbitrary weights, in a future paper~\cite{Plumbings}.
\end{remark}

\begin{proof}
In view of the K\"unneth decomposition for connected sums, see
Theorem~\ref{HolDiskTwo:thm:ConnSumHFm} of~\cite{HolDiskTwo}, it
suffices to consider the case where $G$ is a connected graph.

We will prove inductively that if there is some vertex $v$ in $G$
where $m(v)>d(v)$, then $Y$ is a rational homology sphere and
$\HFa(Y)$ has rank given by the number of elements in $H_1(Y;\Z)$.
(Observe that if this is not the case, and equality holds everywhere,
then it is easy to see by repeated blow-downs that the three-manifold
in question is $S^2\times S^1$, and it is easy to see that
$\HFpRed(S^2\times S^1)=0$, c.f.~\cite{HolDiskTwo}.)

Next, we induct on the number of vertices. Clearly, if the number of
vertices is one, the three-manifold in question is a lens space;
for lens spaces, the conclusion of the theorem
follows easily from the genus one Heegaard diagram 
(c.f. Proposition~\ref{HolDiskOne:prop:Lensspaces} of~\cite{HolDisk}). 

For the inductive step on the number of vertices, we use induction on
$m(v)$ where $v$ is some leaf (vertex with $d(v)=1$). Suppose that
$m(v)=1$. In this case, it is easy to see that $-Y(G)=-Y(G')$, where
$G'$ is the weighted tree obtained from $G$ by deleting the leaf $v$,
and decreasing the weight of the neighbor of $v$ (thought of as a
vertex in $G'$) by one. Observe that $G'$ also satisfies the
hypothesis of the theorem. Thus, the case where $m(v)=1$ follows from
the inductive hypothesis on the number of vertices. More generally,
suppose that $G_1$ is a weighted graph, and we have a leaf $v$ with
$m(v)=k$. In this case, we can form two other weighted graphs $G_2$
and $G_3$, where $G_2$ is obtained from $G_1$
by deleting the leaf $v$, and $G_3$ which
is obtained from $G_1$ by increasing the weight of $v$ by one. 
We have then the
following long exact sequence
(Theorem~\ref{HolDiskTwo:thm:GeneralSurgery} of~\cite{HolDiskTwo}): $$ \begin{CD}
... @>>> \HFa(-Y(G_2))
@>>>\HFa(-Y(G_3)) @>>>
\HFa(-Y(G_1)) @>>> ...  
\end{CD}
$$ By the inductive hypothesis, we know the theorem is true for the
weighted graphs $G_1$ and $G_2$.  Now cobordisms from $-Y(G_2)$ to
$-Y(G_3)$ and from $-Y(G_3)$ and $-Y(G_1)$ (which induce two of the
maps in the above long exact sequence) are clearly
negative-definite. So it follows that $-Y(G_3)$ is a rational homology
sphere, with $$|H_1(Y(G_3);\Z)|=|H_1(Y(G_1);\Z)|+|H_1(Y(G_2);\Z)|.$$
Moreover, by the induction hypothesis, $\HFa(-Y(G_1)$ and
$\HFa(-Y(G_2))$ have no odd-dimensional generators. Since the map
from $\HFa(-Y(G_1))$ to $\HFa(-Y(G_2))$ changes the $\Zmod{2}$
grading, it follows that this map is zero, so that the above long
exact sequence is actually a short exact sequences. This implies that $\HFa(Y(G_3))$
is a free Abelian group with rank
$$\rk \HFa(Y(G_3)) = \rk \HFa(Y(G_1)) + \rk
\HFa(Y(G_2)).$$ The induction hypothesis is equivalent to the
statement that for $i=1,2$, $\HFa(Y(G_i))$ are free and $\rk
\HFa(Y(G_i))=|H_1(Y(G_i);\Z)|$, which in turn gives the corresponding
equation for the graph $G_3$.
\end{proof}

\vskip.2cm
{\noindent{\bf Proof of Theorem~\ref{thm:CertainPlumbings}.}}
According to the definition of $\Phi$, if $X$ is a smooth four-manifold
which can be separated along a rational
homology three-sphere $Y$ into 
$$X=X_1\cup_Y X_2$$
so that $b_2^+(X_i)>0$, then $Y$ constitutes an admissible cut for the definition of 
$\Phi$. If $\HFpRed(Y)=0$, then the invariant $\Phi$ must vanish identically. 
Thus, in this case, the existence of such a decomposition along a graph manifold satisfying
the hypotheses of Theorem~\ref{thm:FloerHomology} gives a vanishing result which
is inconsistent with  Theorem~\ref{intro:NonVanishing}.

In the case where $Y$ is not a rational homology three-sphere, it is formed
as a connected sum of a rational homology three-sphere (as in Theorem~\ref{thm:FloerHomology})
with a collected of copies $S^2\times S^1$. It follows from the behaviour of Floer homology 
under connected sums (c.f.~\cite{HolDiskTwo}) that $\uHFpRed(Y,M)=0$ for any choice of
twisted coefficient system $M$ over $Y$, so we again get a vanishing result for $\Phi$
for any smooth four-manifold which admits the hypothesized decomposition along $Y$.
\qed